\documentclass{article}

\usepackage{arxiv}

\usepackage[utf8]{inputenc} 
\usepackage[T1]{fontenc}    
\usepackage{hyperref}       
\usepackage{url}            
\usepackage{booktabs}       
\usepackage{amsfonts}       
\usepackage{nicefrac}       
\usepackage{microtype}      
\usepackage{graphicx}
\usepackage{natbib}
\usepackage{doi}

\usepackage{amssymb,amsmath,amsthm,latexsym}
\usepackage{natbib}

\usepackage[english]{babel}
\usepackage[utf8]{inputenc}
\usepackage{amsmath}
\usepackage[shortlabels]{enumitem}
\usepackage{graphicx}
\usepackage[colorinlistoftodos]{todonotes}
\usepackage{bbm}  
\usepackage{amsfonts}
\usepackage{lineno}
\usepackage{xspace}
\usepackage{algorithm}
\usepackage{graphicx}

\usepackage{color}
\usepackage{epstopdf}
\usepackage{longtable}
\usepackage{lscape}
\usepackage{indentfirst}
\usepackage{breqn}

\renewcommand{\vec}[1]{\mathbf{#1}}

\newcommand\norm[1]{\left\lVert#1\right\rVert}

\newcommand\Lpnorm[1]{\left\lVert#1\right\rVert _ {L^p(\Omega,\mathbb{R}^d)} }

\newcommand\Lqnorm[1]{\left\lVert#1\right\rVert _ {L^q(\Omega,\mathbb{R}^d)} }

\newcommand\LpOpnorm[1]{\left\lVert#1\right\rVert _ {L^p(\Omega,\mathbb{R}^{d\times d})} }

\newcommand{\expect}[1]{\mathbb{E}{\left[#1\right]}}

\newcommand\LpnormR[1]{\left\lVert#1\right\rVert _ {L^p(\Omega,\mathbb{R})} }

\newcommand\LrnormR[1]{\left\lVert#1\right\rVert _ {L^r(\Omega,\mathbb{R})} }

\newcommand\LphalfnormR[1]{\left\lVert#1\right\rVert _ {L^{p/2}(\Omega,\mathbb{R})} }

\newcommand\Ltwopnorm[1]{\left\lVert#1\right\rVert _ {L^{2p}(\Omega,\mathbb{R}^d)} }
\newcommand\LtwopnormR[1]{\left\lVert#1\right\rVert _ {L^{2p}(\Omega,\mathbb{R})} }

\newcommand \Lfourpnorm[1]{\left\lVert#1\right\rVert _ {L^{4p}(\Omega,\mathbb{R}^d)} }

\newcommand{\Ito}{It{\^o}\xspace}

\newcommand\eval[1]{\mathbb{E}\left[ #1  \right]    }

\newcommand{\innerproduct}[2]{{ \langle} #1 , #2 { \rangle}}
\newcommand{\expint}[1]{\textit{#1}}

\newcommand{\secref}[1]{Section \ref{#1}\xspace}

\newtheorem{definition*}{Definition}
\newtheorem{theorem}{Theorem}

\newtheorem{lemma}{Lemma}
\newtheorem{assumption}{Assumption}

\newcommand{\SG} [1]{\mathbf{\Phi} _{#1} }

\newcommand{\Dt}{\Delta t}

\title{Strong Convergence of a GBM Based Tamed Integrator for SDEs and an Adaptive Implementation}

\author{ \hspace{1mm}{{Utku Erdogan}} \\
Eskisehir Technical University\\
Department of Mathematics\\
Turkey\\
	\texttt{utkuerdogan@eskisehir.edu.tr} \\
	\And
	\hspace{1mm} Gabriel J.  Lord \\
	Department of Mathematics\\
	IMAPP, Radboud University\\
	Nijmegen, The Netherlands\\

	\texttt{gabriel.lord@ru.nl} \\
}

\renewcommand{\shorttitle}{Erdogan and Lord}

\hypersetup{
pdftitle={A template for the arxiv style},
pdfsubject={q-bio.NC, q-bio.QM},
pdfauthor={David S.~Hippocampus, Elias D.~Striatum},
pdfkeywords={First keyword, Second keyword, More},
}

\begin{document}
\maketitle

\begin{abstract}
	 We introduce a tamed exponential time integrator which exploits linear terms in both the drift and diffusion for Stochastic Differential Equations (SDEs) with a one sided globally Lipschitz drift term. Strong convergence of the proposed scheme is proved, exploiting the boundedness of the geometric Brownian motion (GBM) and we establish order 1 convergence for linear diffusion terms. In our implementation we illustrate the efficiency of the proposed scheme compared to existing fixed step methods and utilize it
  in an adaptive time stepping scheme. Furthermore we extend the method to nonlinear diffusion terms and show it remains competitive.
  The efficiency of these GBM based approaches are illustrated by considering some well-known SDE models.
\end{abstract}

\keywords{First keyword \and Second keyword \and More}

\section{Introduction}
We consider the numerical integration of Stochastic Differential Equations (SDEs) 
\begin{equation}
\label{eq:genSDE}
    dX_t=\mu(X_t)dt +\sigma(X_t)dW_t, \quad X_0=\xi 
\end{equation}
on finite time interval $[0,T]$ where $\mu:\mathbb{R} ^d \rightarrow \mathbb{R} ^d$ and  $\sigma:\mathbb{R} ^d \rightarrow \mathbb{R} ^{d \times m}$  and $W_t$  is  m dimensional Wiener process   on probability space  $(\Omega , \mathcal{F},\mathbb{P})$   with filtration $({\mathcal{F} _ {t}  }) _{t \in 
  [0,T]}$. This has a well established theory (see for example \cite{kloeden2011}) for globally Lipschitz drift and diffusion coefficients $\mu$  and $\sigma$.  
  We propose and prove convergence of a tamed method for one-sided Lipschitz drift terms $\mu$ with linear drift and diffusion terms, see \eqref{eq:EqAFB} and extend the method to equations of the general form \eqref{eq:genSDE} in Section \ref{sec:LargerClass}. 
     
There has been great interest recently in methods for SDEs with one-sided Lipschitz drift and we give a brief review of alternative methods and approaches. First we note that when a numerical method is drift implicit then the nonlinearity $\mu$ can be controlled, see for example \cite{HighamMaoSzpruch2013,YaoGan18}, but at the expense of a linear solve at each step. 
 Higham et al.
 \cite{higham2002} investigated 
  the strong convergence of explicit Euler-Maruyama (EM) scheme under non-globally Lipschitz coefficients and concluded that if the numerical solution has bounded moments then strong convergence can be obtained. However, Hutzenthaler et al.  \cite{hutzenthaler2011} proved  that the explicit EM scheme produces numerical solution with  unbounded moments for  superlinearly growing coefficients $\mu$ and $\sigma$. 
  Additionally in \cite{hutzenthaler2012} they proposed an explicit scheme, "tamed EM", based on replacing a superlinearly growing coefficient with its first order bounded approximation at each step in the scheme. For this scheme they then proved moment bounds and strong convergence for  SDEs with one sided globally Lipschitz  drift coefficient $\mu$. There follows a number of different taming type schemes such as  \cite{sabanis2013note,sabanis2016,MaoStopped2013,Izgi2018,zong2014}.
  We mention in particular Wang and Gan \cite{wang} who devised a tamed version of Milstein for SDEs with commutative noise that we consider later Section \ref{Sec:Num} and the generalization considered in the work of \cite{KumarSabanis}.

Another approach to deal with superlinearly growing coefficients is to  modify the numerical solution before each iteration. Mao et al, see for example \cite{TruncatedMao2018,TruncatedMao2015}, devised truncated methods obtaining numerical solution by truncating intermediate terms according to growth rate of $\mu$  and $\sigma$  to avoid blow up.  
In this category, Beyn et al. \cite{Beyn2016,Beyn2017} introduced the Projected Euler and Milstein methods consisting of the classical EM/Milstein methods with a projection of the solution onto a ball whose radius is proportional to a negative power of step size.  To prove strong convergence of proposed schemes, they utilized suitable generalization of the concepts C-stability and B-consistency which require less strict conditions on moments of numerical solutions.

In  \cite{kelly2016adaptive},  Kelly and Lord proposed adaptive EM methods, extended to Milstein in \cite{KellyMilsteinArxiv}, which aim to control the growth of the numerical solution by choosing adaptively a suitable time step. After a specified cutoff a 'backstop' fixed step method is used.
We later examine in Section \ref{Sec:Num} the numerical performance of adaptive methods based on the schemes here and \cite{utkuLord}.

We propose novel explicit tamed based methods that are in the broad class of exponential integrators 
and build on the work in \cite{utkuLord}. This takes advantage of linear terms not only in the drift but also in the diffusion 
which is not the case in earlier exponential integrators such as \cite{biscay1996local,Lord2004587}, see also the short review in \cite{utkuLord}.
  Recent related work on exponential integrators for SDEs includes \cite{doi:10.1137/15M1041341} looking at second-order weak convergence for Runge-Kutta methods,  \cite{YangMaLi} who propose Magnus type exponential integrators for Stratonovich SDEs and \cite{DebrabantKvaernoMattsson} that examines families of Runge-Kutta Lawson schemes and their weak and strong convergence in the absence of a nonlinear drift term.
Specifically we examine new explicit tamed based methods for the semilinear SDE  
\begin{equation} \label{eq:EqAFB}
  dX_t= \left( A X_t+F(X_t) \right) dt + \sum _{i=1}
  ^m  \left(  B_i X_t\right)  dW_t^i
  \quad X_0=\xi \in  \mathbb{R} ^d, 
\end{equation}
with $m$, $d\in\mathbb{N}$ 
where $W_t^i$ are iid Brownian motions, $F :
\mathbb{R} ^d \rightarrow \mathbb{R} ^d$  is one sided 
Lipschitz (see Assumption \ref{ass:1}), and the matrices $A$, $B_i \in
\mathbb{R} ^{d \times d }$  satisfy the following zero commutator conditions 
\begin{equation} \label{eq:commute}
AB_i-B_iA=0, \quad B_jB_i-B_iB_j=0
\qquad \text{for} \quad  i,j=1 \hdots m.
\end{equation}
When these commutator conditions hold
then we can exploit the exact flow of \eqref{eq:EqAFB}.
SDE's of this form arise in many areas of  engineering and science, 
and in particular an SDE with diagonal noise satisfies the commutator condition.
The standard cubic nonlinear scalar SDE, see \eqref{eq:GL}, is a classic example but we also consider in Section \ref{Sec:Num} examples from the mathematical biology: stochastic Lotka--Volterra and HIV models as well as tumor growth model. In addition
the SDEs that arise from the spectral discretization of semilinear Stochastic Partial Differential Equations (SPDEs) often fall into this class.

The method we examine is based on the exact flow of \eqref{eq:EqAFB} which is expressed in the following integral form
\begin{equation}  \label{eq:exact}
X_{t_{n+1}}=\SG{t_{n+1},t_n} \left( X_{t_{n}} + \int _{t_n} ^{t_{n+1}} \SG {s,t_n} ^{-1} F \left( X_s \right) ds \right),
\end{equation}
where the fundamental matrix 
\begin{equation} \label{eq:semigroup}
\SG{t,t_0}=\exp \left(( A-\frac{1}{2} \sum_{i=1} ^m B_i ^2)(t-t_0) +
  \sum_{i=1}^m B_i (W_t^i  -W_{t_0}^i) \right)
\end{equation}
is the solution to the linear equation 
\begin{equation} \label{eq:Homogen}
d  \SG{t, t_0} = A  \SG{t,t_0}   dt + \sum _{i=1} ^m  \ B_i \SG{t,t_0} dW_t^i, \qquad \SG{t_0, t_0} =I_d.
\end{equation} 
In \cite{utkuLord} we proposed the numerical scheme \expint{EI0}
\begin{equation}  \label{eq:EI0}
Y_{n+1}^N=\SG{t_{n+1},t_n}\left( Y_n^N+ F(Y_n^N) \Delta t  \right),
\end{equation}
where $\Delta t=\frac{T}{N}$ defined by the uniform time partition $0=t_0<t_1    <t_2<\hdots <t_N=T$.
Here we introduce a tamed version of \eqref{eq:EI0} that we denote \expint{TamedEI0}
\begin{equation}  \label{eq:tamedEI0}
    Y_{n+1}^N=\SG{t_{n+1},t_n}\left( Y_n^N+ \tilde{F}
      (Y_n^N) \Delta t  \right), \qquad Y_0^N=\xi
  \end{equation}
where 
\begin{equation}
\label{eq:taming}
    \tilde{F}(Y_n^N)=\frac{F(Y_n^N)}{1+\Delta t \norm{F(Y_n^N)}}.
\end{equation}
Compared to method in \cite{utkuLord} the taming has changed the scheme and a key step is to obtain moment bounds on the numerical method.
It is clear from \eqref{eq:tamedEI0} that the Brownian increment term is included in the $\SG{t_{n+1},t_n}$. As a consequence the  convergence analysis needs to be adapted from that in \cite{hutzenthaler2012,wang,zong2014,sabanis2013note} and we need to exploit the boundedness of geometric Brownian motion (GBM). 
The numerical results of \cite{utkuLord} indicate that \expint{EI0} was the most efficient of the three schemes considered. This remains the case for taming versions of these and hence we only consider \expint{EI0} here.
We restrict our analysis to the linear diffusion case, as this leads to an improved rate of convergence, however we examine nonlinear diffusion numerically.
Our main result in the paper is a strong convergence of the scheme that we state and prove in Section \ref{sec:StrongConv} along with some necessary lemmas.  In Section \ref{sec:AssumptionPreliminary} we introduce our notation and give a preliminary result on the linear system before in Section \ref{sec:boundedmoment} stating results for the boundedness of the moments of the scheme (the proofs are given in \ref{app:boundedmoment}).
In Section \ref{sec:backstop} we discuss how the scheme \eqref{eq:tamedEI0} can be used as the backstop method for an adaptive time stepping strategy and in 
Section \ref{sec:LargerClass} we discuss how \eqref{eq:tamedEI0} can be extended to general systems of the form \eqref{eq:genSDE}. Finally in Section \ref{Sec:Num} we illustrate convergence and compare the efficiency of methods. We observe that amongst the fixed step methods \eqref{eq:tamedEI0} is the most efficient on the examples with linear diffusion and that the most efficient is an adaptive version. 
Where the diffusion is no longer linear our proposed schemes remain competitive.

\section{Assumptions and preliminary results}
\label{sec:AssumptionPreliminary}
 
We first introduce the notation for the norms we work with. We let $\norm{.}$ denote the standard Euclidean norm for vectors in $\mathbb{R}^d$ and also for the subordinate matrix norm for matrices in $\mathbb{R}^{d\times d}$. In addition we let   $\norm{.}_{T_3}$ denote the subordinate tensor norm  of rank 3 in   $\mathbb{R}^{d\times d \times d}$
given by
$$
\|J(x)\|_{T_3} = \sup_{\|h_1\|\leq 1,\|h_2\|\leq 1}
\|J(x)(h_1,h_2)\|.$$
Finally we introduce the norm $\Lpnorm{.}$ where 
$\Lpnorm{ v}:= \left(\expect{\norm{v}^p}\right)^{1/p}$.
Throughout the paper, we make use of the following time regularity and norm properties of GBM.
\begin{lemma}\label{lem:GBM_half_order}
Let $0<t-s<1$ and consider  the stochastic matrix  $\SG {t,t_0}$  given in \eqref{eq:semigroup} with $t_0=s$.
Then for each ${\mathcal{F} _ {s}}$ measurable random variable $v \in$  $L^p(\Omega,\mathbb{R}^d)$  with  $ p\geq 2 $
\begin{align} \label{eq:halfOrder}
 \Lpnorm{\SG {t,s}v-v} \leq C_p \vert t-s \vert ^{1/2},\\ \label{eq:halfOrder_inverse}
\Lpnorm{\SG {t,s} ^{-1} v-v} \leq C_p \vert t-s \vert ^{1/2}.
\end{align}
\end{lemma}

\begin{proof}
For \eqref{eq:halfOrder}, by definition of GBM,  $U_t=\SG{t,s}U_s$ solves 
\begin{equation} \label{eq:GBM_IntegralEq}
  U_t=U_s+ \int_{s} ^t A U_r  dr +\int_{s} ^ t \sum _{i=1}
  ^m  B_i U_r dW_r^i .
\end{equation}
Taking the pth power of the Euclidean norm 
we have
\begin{equation*}
 \norm{U_t-U_s}^p \leq  2^{p-1} \norm{\int_{s} ^t A U_r dr}^p + 2^{p-1}\norm{\int_{s} ^ t \sum _{i=1}
  ^m  B_i U_r dW_r^i }^p.
\end{equation*}
By taking the expectation, using the subordinate matrix norm inequality, Jensen's inequality and the moment inequality in \cite[Theorem 1.7.1]{maoBook} we find
\begin{multline*}
 \Lpnorm{U_t-U_s}^p \leq 2^{p-1} (t-s)^{p-1} \int_{s} ^t \norm{A}\Lpnorm{U_r }^p dr \nonumber \\+2^{p-1} C(p,m)(t-s)^{p/2-1} \int_{s} ^ t \sum _{i=1}
  ^m  \norm{B_i}\Lpnorm{ U_r}^p dr.
\end{multline*}
By boundedness of $U_t$, the solution to   \eqref{eq:GBM_IntegralEq}, in $L^p (\Omega,\mathbb{R}^d )$, we obtain \eqref{eq:halfOrder}.

To prove \eqref{eq:halfOrder_inverse} we note that the  operator $U_t=\SG{t,s}^{-1}U_s$ satisfies the linear  SDE
\begin{equation}
dU_t=(-A+\sum _{i=1} ^m  B_i^2) U_tdt-\sum _{i=1}
  ^m  B_i U_t dW_t^i 
\end{equation}
and the desired result follows as for \eqref{eq:halfOrder}.
 \end{proof}

We now prove a preliminary lemma that shows the subordinate matrix norm  $\LpOpnorm { \SG{t,s}}$  of the operator $\SG {t,s}$ acting on  the space $L^{p}(\Omega, \mathbb{R}^d)$   is well defined and bounded on a finite time interval.

\begin{lemma} \label{lem:SGBoundedness}
Suppose $p \geq 2$, $s<t$ and let $\SG {t,s}$ be given in \eqref{eq:semigroup} with $t_0=s$. 
\begin{itemize}
    \item[i)]

For any $\mathcal{F}_{s}$ measurable random variable $v$ in  $L^p (\Omega,\mathbb{R}^d )$
  $$ \Lpnorm{\SG {t,s} v}  \leq    \exp{ \left(\left(\norm{A}+\frac{p-1}{2}\sum _{i=1}
  ^m  \norm{ B_i}^2\right)   (t-s)\right)  }  \Lpnorm{ v}.$$ 
  \item[ii)]
   For any $\mathcal{F}_{t}$ measurable random variable $v$ in $L^p (\Omega,\mathbb{R}^d )$

  $$ \Lpnorm{\SG {t,s} v} \leq K \Lqnorm{ v}, $$ 
  where $K_p=\LrnormR{\norm{\SG{t,s}}} < \infty $ for  $1/p=1/r+1/q$ with $r$, $q>1$.
\end{itemize}
\end{lemma}
\label{app:AssumptionPreliminary}
\begin{proof}
\begin{itemize}
\item[i) ]By definition 
$U_t=\mathbf{\Phi}_{t,s} v$ solves \eqref{eq:GBM_IntegralEq} 
where $v$ is an ${\mathcal{F} _ {s}}$ measurable random variable in $L^p(\Omega,R^d)$.
By \Ito's lemma,
\begin{align}
 \norm{U_t}^p=& \norm{v}^p+\int_{s} ^{t} p \norm{U_r}^{p-2} \innerproduct{U_r}{AU_r}dr  \nonumber \\ 
 & +\frac{1}{2} \int_{s} ^{t} \left( p \norm{U_r}^{p-2}\sum _{i=1}
  ^m  \norm{ B_i U_r}^2 +p(p-2)\norm{U_r}^{p-4}\sum _{i=1}
  ^m \vert \innerproduct{B_i U_r}{U_r} \vert ^2 \right)dr \nonumber \\
   &  + \int_{s} ^{t} p \norm{U_r}^{p-2}\sum _{i=1}
  ^m  \innerproduct{B_i U_r}{U_r} dW_r^i.
\end{align}
By applying the Cauchy-Schwarz inequality, taking the expectation and using mean zero property of \Ito integrals we have
\begin{align*}
\eval{ \norm{U_t}^p}\leq & \eval{\norm{v}^p}+\left(p\norm{A}+\frac{p}{2}\sum _{i=1}
  ^m  \norm{ B_i}^2 +\frac{p(p-2)}{2} \sum _{i=1}
  ^m  \norm{ B_i}^2\right)\int_{s} ^{t} \eval{ \norm{U_r}^{p}} dr.  
\end{align*}
An application of the standard Gronwall inequality, gives
$$
\eval{ \norm{U_t}^p} \leq   \exp\left({ p\left(\norm{A}+\frac{p-1}{2}\sum _{i=1}
  ^m  \norm{ B_i}^2 \right)  (t-s) }\right)    \eval{ \norm{v}^p}
$$
and hence the result.

\item[ii)] By applying the subordinate  norm inequality to $\norm{\SG {t,s} v}$, taking the expected value of both sides, and finally  applying H\"{o}lder's inequality in  $L^p (\Omega,\mathbb{R}^d )$, we have the inequality
\begin{align}
\eval{\norm{\SG {t,s} v} ^p } ^{1/p} & \leq  \eval{\norm{\SG {t,s}} ^p    \norm{ v} ^p  } ^{1/p} \\
           & \leq   \eval{\norm{\SG {t,s}} ^r }^{1/r}  \eval{\norm{v} ^q }^{1/q}.
\end{align}
It remains to show that $\eval{\norm{\SG {t,s}} ^r }< \infty$. Note that
\begin{equation}
\norm{\SG {t,s}} ^r \leq  \exp\left( \norm{ r( A-\frac{1}{2} \sum_{i=1} ^m B_i ^2)(t-s) } \right) \exp \left( 
  r  \sum_{i=1}^m \norm{ B_i} \vert W_t^i  -W_{s}^i  \vert  \right). 
\end{equation}
By taking  expected value and considering the  independence of  random  values $\vert W_t^i  -W_{s}^i  \vert$, for $i=1,2,\hdots m$, we have
\begin{eqnarray*}
\lefteqn{\LrnormR{\norm{\SG {t,s}}} ^r } \\
& \leq & \exp\left( \norm{ r( A-\frac{1}{2} \sum_{i=1} ^m B_i ^2)(t-s) } \right)  \prod_{i=1}^m  \eval{ \exp \left( 
  r  \norm{ B_i} \vert W_t^i  -W_{s}^i  \vert  \right) } \\
&  = &\exp\left( \norm{ r( A-\frac{1}{2} \sum_{i=1} ^m B_i ^2)(t-s) } \right)  \prod_{i=1}^m  \exp(r^2 \norm{B_i}^2 (t-s)) \left(1+\text{erf}(\frac{r \norm{B_i}}{\sqrt {2} }) \right)\\
& < & \infty
\end{eqnarray*}
\end{itemize}
\end{proof}

\begin{assumption}\label{ass:1} Let  $F$ and its Jacobian  $\vec{D}F$ be continuously differentiable. Furthermore let there exist positive constants  $K\geq 1$  and $c$, $q\geq 1$ such that  $\forall x,y \in \mathbb{R}^d$
\begin{enumerate}[i)]
\item $ \innerproduct{x-y}{ F(x)-F(y) }\leq K \norm{x-y}^2 $
\item $\norm{ \vec{D}F(x) } \leq  K (1+\norm{x}^c )$
\item $\norm{\vec{D}^2 F(x) }_{T_3} \leq  K (1+\norm{x}^q )$.
\end{enumerate}  
\end{assumption}

We note that \cite{Beyn2017} and \cite{KumarSabanis} both examine the case where $F$ is only continuously differentiable and there is no condition on $D^2F$. We use this condition as in \cite{wang} to get a martingale to apply Doob's maximal inequality. It may be possible to weaken this assumption on $F$ but would require an alternative approach to the proof. Our proof is for the case of a linear diffusion, however in Section \ref{sec:LargerClass} we extend to the nonlinear case and we comment on the expected rate of convergence there.

Our focus in the analysis is on obtaining an order $1$ rate of convergence and hence we analyse the special diffusion term in \eqref{eq:EqAFB} and only consider more general diffusion terms (such as in \cite{wang,Beyn2017,KumarSabanis}) in our numerical experiments.

\section{Boundedness of moments of the numerical scheme}
\label{sec:boundedmoment}

The main idea of the proof of the moment bounds is adapted from that in \cite{hutzenthaler2012,wang,zong2014,sabanis2013note}. We highlight here (and in the proofs in \ref{app:boundedmoment}) the differences required to examine \eqref{eq:tamedEI0}. In the standard way we start by establishing the boundedness of numerical solutions on a set $\Omega_n^N$ where increments of the noise are controlled. 
Then, by using standard inequalities on the measure of complement set, we obtain the boundedness of $p$th moments of numerical solution. We make extensive use of Lemma \ref{lem:SGBoundedness} and the key property of the taming factor that
\begin{equation}\label{eq:tamebound}
\norm{\tilde{F}(Y_n^N)} = \dfrac{\norm{F(Y_n^N)}}{1+\Delta t \norm{F(Y_n^N)}}\leq \norm{F(Y_n^N)}.
\end{equation}
We introduce appropriate sub events of $\Omega$. 
We let $\Omega_0^N=\Omega$, then for $n\in\{1,2,\ldots, N\}$
\begin{equation}\label{eq:Omega}
\Omega_n^N=\left\{\omega \in \Omega \vert ,\sup _{0\leq k \leq n-1} D_k^N(\omega) \leq N^{1/2c}
,\sup _{0\leq k \leq n-1} \norm{\sum_{i=1}^m B_i\Delta W_{t_k}^i}\leq 1 \right\}, 
\end{equation}
 and introduce the dominating stochastic process,  $D_0^N=\left( \lambda + \norm{\xi} \right) e^{\lambda} $ and 
\begin{equation}\label{eq:Dn}
D_n^N=\left( \lambda + \norm{\xi} \right) e^{\lambda} \sup_{0\leq u \leq n} \prod_{k=u}^{n-1}\norm{\SG {t_{k+1},t_k}}, \qquad n \in  \lbrace 1,2,\hdots, N \rbrace
\end{equation}
where
$$\lambda=e^{1+T \norm{ A-\frac{1}{2} \sum_{i=1} ^m B_i ^2}} \left(1+4TK+2T \norm{F(0)} +2K \right) ^2. $$
The first result shows we can dominate the numerical solution on the set $\Omega_n^N$.
\begin{lemma} \label{lem:boundedness1}
Let $Y_n^N$ be given by \eqref{eq:tamedEI0}.
For all $n=0,1,\hdots, N$ we have $$1_{\Omega_n^N} \norm{Y_n^N} \leq D_n^N.$$%
\end{lemma}
The second result then bounds the dominating process.
\begin{lemma} \label{lem:Dboundedness} For all $p \geq 1$,
$
\sup_{N \in \mathbb{N} } \eval{\sup_{0\leq n\leq N} \vert D_n^N \vert ^p} < \infty .
$
\end{lemma}

For the complement we have from \cite[Lemma 3.6]{hutzenthaler2012} the following:
\begin{lemma} \label{lem:complementProbability}
For all $p\geq1$, $
\sup_{N \in \mathbb{N}} (N^p \mathbb{P} [(\Omega _N^N)^c] ) < \infty.
$
\end{lemma}
Finally we obtain bounded moments for the numerical scheme \eqref{eq:tamedEI0}. 
\begin{theorem}\label{teo:boundedMoment}
Let $Y_n^N$ be given by \eqref{eq:tamedEI0}.
Then, for  all  $p \in [1,\infty)$
$$
\sup_{N \in \mathbb{N}} \left[\sup_{0\leq n\leq N}  \eval{\norm{Y_n^N}^p} \right]<\infty.
$$
\end{theorem}
The proofs of Theorem \ref{teo:boundedMoment} and relevant lemmas are contained in  \ref{app:boundedmoment}.

\section{Strong convergence of TamedEI0}
\label{sec:StrongConv}
The aim of this section is to prove the following strong convergence result
 $$
\left( \eval{\sup_{t \in [0,T]} \norm{X_t-\bar{Y}_t}^p}  \right)^{\frac{1}{p}} \leq  C_{p,T} \Delta t,
$$
 where $C_{p,T}$ is a constant independent of $\Delta t$, see Theorem \ref{thm:convergence} in \secref{sec:main}.
Before we state and prove this result we require a number of preliminary lemmas.

\subsection{Preliminary Lemmas}\label{sec:prelim_main}

\begin{lemma} \label{lem:boundedFexact}
Let Assumption \ref {ass:1} hold. Then for all $p\geq 1$, we have
$$
\sup_{t \in [0,T]}  \Lpnorm{X_t} < \infty, \qquad
\sup_{t \in [0,T]}  \Lpnorm{F(X_t)} < \infty.  
$$
\end{lemma}
\begin{proof}
These follow from the boundedness of moments, 
see 
\cite[Theorem 4.1]{maoBook} and \cite[Lemma 3.4]{wang}.
\end{proof}
Now we define a continuous form $\bar{Y}_t$ of the numerical solution by introducing  
\begin{equation} \label{eq:hatnotation}
\hat{t}=t_n, \text{ for } t_n\leq t< t_{n+1},
\end{equation}
and setting  
 \begin{equation} \label{eq:final_cont_form}
 \bar{Y}_t= \SG{t,0} \xi+\SG{t,0}  \int_0^t \SG{\hat{s},0}^{-1}  \tilde{F}( \bar{Y}_{\hat{s}})ds.
 \end{equation}
Note that $\bar{Y}_t$ agrees with the approximation $Y_n^N$ at  $t=t_n$.

The following lemma provides an integral representation for the continuous version of numerical solution.
 \begin{lemma}\label{lem:itoEq}
Let $ \bar{Y}_t$ be the interpolated  numerical  solution given in \eqref{eq:final_cont_form}. Then the  differential for this solution is given by
\begin{equation}
\label{eq:cts}
d\bar{Y}_t=\left(     A  \bar{Y}_t +\SG{t,\hat{t} } \tilde{F} (\bar{Y}_t )          \right) dt+\sum_{i=1} ^{m}B_i  \bar{Y}_t dW_t^i.
\end{equation}
\end{lemma}
\begin{proof}
By definition of GBM
$$
d \SG{t,0} ^{-1}=(-A+\sum _{i=1}
  ^m B_i^2) \SG{t,0} ^{-1} dt - \sum _{i=1}
  ^m B_i \SG{t,0} ^{-1} dW_t^i.
$$ 
Let us set
$d \bar{Y}_t=\mu dt+\sum _{i=1}
  ^m\sigma_i dW_t^i$
and seek the appropriate $\mu$ and $\sigma_i$.
The  representation for matrix-vector product gives us
\begin{equation}
\label{eq:tmp1}
d\left[  \SG{t,0} ^{-1} \bar{Y}_t \right]= \SG{t,0} ^{-1}\left[(-A+\sum _{i=1}
  ^m B_i^2) \bar{Y}_t+\mu-\sum _{i=1}
  ^m B_i\sigma_i \right] dt+ \SG{t,0} ^{-1}\left[-\sum _{i=1}
  ^mB_i \bar{Y}_t +\sigma_i \right] dW_t^i. 
\end{equation}
On the other hand, the linearly interpolated continuous solution is given by
$$
\bar{Y}_t = \SG{t,0} \xi+  \SG{t,0}\int_{0} ^ {t}  \SG{\hat{s},0}^{-1} \tilde{F} (\bar{Y}_{\hat{s}})  ds.
$$
By comparison of the above expression with \eqref{eq:tmp1} after multiplying through by $\SG{t,0}$, we find 
$$
   \SG{t,0} ^{-1}\left[(-A+\sum _{i=1}
  ^m B_i^2) \bar{Y}_t+\mu-\sum _{i=1}
  ^m B_i\sigma_i \right]  = \SG{\hat{t},0}^{-1}  \tilde{F}(\bar{Y}_{\hat{t}})
$$
and 
$$
\SG{t,0} ^{-1}\left[-\sum _{i=1}
  ^mB_i \bar{Y}_t +\sigma_i \right] =0.
$$
Finally, we have that
$
\mu=A \bar{Y}_t + \SG{t,\hat{t}} \tilde{F}(\bar{Y}_{\hat{t}})$
and 
$\sigma_i=B_i \bar{Y}_t.$
\end{proof}

\begin{lemma} \label{lem:boundedFnumeric}
Let Assumption \ref {ass:1} hold. Then for all $p\geq 1$, we have the following: 
\begin{enumerate}[i)]
\item For $Y_n^N$
$$
\sup_{N \in \mathbb{N}} \sup_{0\leq n \leq N}   \eval{\norm{F(Y_n^N)^p} } < \infty$$ 
$$\sup_{N \in \mathbb{N}} \sup_{0\leq n \leq N} \  \eval{ \norm{\vec{D}F( Y_n^N)}^p } < \infty,\quad \text{and} \quad \eval{ \norm{\vec{D}^2F( Y_n^N)}_{T_3}^p } < \infty. 
$$
\item For $\bar{Y}_t$ 
$$
\sup_{t \in [0,T]}  \Lpnorm{\bar{Y}_{t}} < \infty, \quad \text{and} \quad\sup_{t \in [0,T]}  \Lpnorm{F(\bar{Y}_{t})} < \infty. $$
\end{enumerate}
\end{lemma}
\begin{proof}
Theorem \ref{teo:boundedMoment} and the polynomial growth condition on  $\vec{D}F$ and $ \vec{D}^2 F$ from Assumption \ref{ass:1} imply the  estimates in i). 
 
For the estimates in ii) we  consider the  one step continuous  extension
\begin{equation} \label{eq:OneStepContExt}
 \bar{Y}_{t}=\SG{t,\hat{t}} \left( \bar{Y}_{\hat{t}}+\int_{\hat{t}}  ^t  \tilde{F}( \bar{Y}_{\hat{t}} ) d s \right)
\end{equation}
 where $\hat{t}$ is as defined in \eqref{eq:hatnotation}.
 Therefore,  by Theorem \ref{teo:boundedMoment}, we see $\Lpnorm{\bar{Y}_{\hat{t}}}=\Lpnorm{Y _n^ N }$  and for $p\geq 2$
$$\Lpnorm{ \bar{Y}_{t}} \leq  e^{ \left(\norm{A}+\frac{p-1}{2}\sum _{i=1}
  ^m  \norm{ B_i}^2 \right)  (t - {t_n}) }   \left(\Lpnorm{ \bar{Y}_{\hat{t}} }+\int_{\hat{t}}  ^t  \Lpnorm{\tilde{F}( \bar{Y}_{\hat{t}} )} d s \right)<\infty.
 $$
 By H\"older's inequality we obtain this for $p\geq 1$.
Similarly, the polynomial growth condition on  $\vec{D}F$ and the Mean Value Theorem implies the final estimate.
\end{proof}

\begin{lemma} \label{lem:orderOfIncrementsInLp}
Let  Assumption \ref {ass:1} hold on $F$. Then for all $p\geq 2 $, $t-\hat{t}<1$,
$$ \Lpnorm{\bar{Y}_t - \bar{Y}_{\hat{t}}} \leq C_p (t-\hat{t}) ^{1/2} $$
and 
$$\Lpnorm{F(X_t) - F(X_{\hat{t}}) } \leq C_p (t-\hat{t}) ^{1/2}.$$
\end{lemma}

\begin{proof}

The time difference in  one step  continuous form of the numerical solution can be written as
 \begin{equation}\label{eq:timeDifference}
 \bar{Y}_t-\bar{Y}_{\hat{t}}= \left(\SG{t,\hat{t}}-I\right)\bar{Y}_{\hat{t}} +  (t-\hat{t}) \SG{t,\hat{t}}  \tilde{F}( \bar{Y}_{\hat{t}})   
 \end{equation}
since
$$
 \bar{Y}_t=\SG{t,\hat{t}} \left( \bar{Y}_{\hat{t}}+\int_{\hat{t}}  ^t  \tilde{F}( \bar{Y}_{\hat{t}} ) d r \right).
$$
Taking  norm  on  the space $L^p (\Omega,\mathbb{R}^d )$, we have, 
 \begin{equation}
 \Lpnorm{\bar{Y}_{t}-\bar{Y}_{\hat{t}}} \leq \Lpnorm{\SG{t,\hat{t}}\bar{Y}_{\hat{t}}-\bar{Y}_{\hat{t}} }+  (t-\hat{t}) \Lpnorm{\SG{t,\hat{t}}  \tilde{F}( \bar{Y}_{\hat{t}}) }.  
 \end{equation}
 Using \eqref{eq:halfOrder} for the first term above we obtain the first estimate of the lemma. For the second estimate on the exact solution, see  \cite[Lemma 3.5]{wang}. 
\end{proof}

\begin{lemma} \label{lem:GBM_F_order}  
Let Assumption \ref{ass:1} hold on F. Consider the GBM operator $\SG{t,\hat{t}}$ from \eqref{eq:semigroup} (with $t_0=\hat{t}$) and its inverse  $\SG{t,\hat{t}} ^{-1}$ where $t-\hat{t}<1$. Then for $ p\geq 2 $

\begin{itemize}

\item [i)] $\Lpnorm{\SG{t,\hat{t} } F (X_t) - F(X_t)} \leq C_p (t-\hat{t} )^ {1/2} $ 

\item[ii)]  $\Lpnorm{\SG{t,\hat{t} } ^{-1} F (X_t) - F(X_t)} \leq C_p (t-\hat{t} )^ {1/2}$.
\end{itemize} 
\end{lemma}
\begin{proof}
Since $X_t$  is  not  $\mathcal{F}_{\hat{t}}$ measurable, we are not able to apply Lemma \ref{lem:GBM_F_order} directly. By adding and subtracting $F( X_{\hat{t}} )$, we  have
\begin{multline}\Lpnorm{\SG{t,\hat{t} } F (X_t) - F(X_t)}  \leq     \Lpnorm{\SG{t,\hat{t} } F (X_t) - F(X_{\hat{t}})} \\ + \Lpnorm{  F(X_{\hat{t}}) - F(X_t)  }.
\end{multline}
  The second term on RHS above is bounded as desired due to Lemma \ref{lem:orderOfIncrementsInLp}. Let us consider the first term. By a Taylor expansion of $F$ about $X_{\hat{t}}$  and the linearity of $\SG{t,\hat{t} }$, we have
\begin{multline*} 
\Lpnorm{\SG{t,\hat{t} } F (X_t) - F(X_{\hat{t}})} \leq \Lpnorm{\SG{t,\hat{t} } F (X_{\hat{t}}) - F(X_{\hat{t}})} \\+  \Lpnorm{\int_0^1 \vec{D}F(Z_{\theta}) \left(\bar{X}_t-\bar{X}_{\hat{t}}\right) d\theta }
\end{multline*}
where $ Z_{\theta}:= X_{\hat{t}}+\theta \left(X_t-X_{\hat{t}}\right)$, $\theta\in[0,1]$.
The first term provides the desired order by \eqref{eq:halfOrder} of Lemma \ref{lem:GBM_half_order}. For the second term, we have
\begin{eqnarray}
\lefteqn{\int_0^1 \Lpnorm{  \vec{D}F(Z_{\theta})  \left(X_t-X_{\hat{t}}\right) } d \theta } \nonumber \\
    &\leq&  \int_0^1\LpnormR{ \norm{  D F\left(Z \right) }  \norm{ X_t-X_{\hat{t}} } } d \theta  \nonumber  \\
        &  \leq & K \LpnormR{ \Big(  1+\norm{X_{\hat{t}}}^q +\norm{X_t}^q   \Big)    \norm{ X_t-X_{\hat{t}}  } }  \nonumber \\
 &  \leq & K \LtwopnormR{   1+\norm{X_{\hat{t}}}^q +\norm{X_t}^q   }   \LtwopnormR{ X_t-X_{\hat{t}}   }  \nonumber \\
 & \leq & C \sqrt{t-\hat{t}} ,
\end{eqnarray}
by bounded moments and the time regularity  of the exact solution to the SDE \eqref{eq:exact}, \cite[Theorem 4.1]{maoBook}. 
\end{proof}

Now, we are going to consider the order of three remainder terms arising from \Ito-Taylor expansions. 
First we introduce the remainder term $R_{GBM}(t,\hat{t})$ from the \Ito-Taylor expansion of GBM near $t=\hat{t}$:  

\begin{equation}\label{eq:R_GBM}
 \left( \SG{t,\hat{t}}-I\right)\bar{Y}_{\hat{t}} =\sum _{i=1}^m   B_i   \bar{Y}_{\hat{t}}\left(  W_t^i-W^i_{\hat{t}}\right)+R_{GBM}(t,\hat{t}).
\end{equation}
\begin{lemma}\label{lem:RemainderGBM}
Let Assumption \ref{ass:1} hold an $t-\hat{t}<1$. Then for $ p\geq 2 $,
$$
\Lpnorm{R_{GBM}(t,\hat{t})}   <C_p (t-\hat{t}).
$$
\end{lemma}

\begin{proof}
By the \Ito-Taylor expansion of GBM, we have
\begin{multline} \label{eq:remainderGBM}
R_{GBM}(t,\hat{t})=\sum _{i,j=1} ^m    B_iB_j  \int _{\hat{t}} ^{t} \int _{\hat{t}} ^{s} \SG{z,\hat{t}} \bar{Y}_{\hat{t}}dW_i(z)dW_j(s) +  \int _{\hat{t}} ^{t} A \SG{s,\hat{t}} \bar{Y}_{\hat{t}}ds \\
+\sum _{i=1} ^m    B_iA  \int _{\hat{t}} ^{t} \int _{\hat{t}} ^{s} \SG{z,\hat{t}} \bar{Y}_{\hat{t}}dzdW_s^i.
\end{multline}
Applying the moment inequality in  \cite[Theorem 1.7.1]{maoBook}, we obtain the result. 
\end{proof}
The second remainder term comes from deterministic Taylor expansion of $ F(\bar{Y}_t)$ near $F( \bar{Y}_{\hat{t}})$: 
\begin{equation} \label{eq:deterministicTaylor}
F(\bar{Y}_t)- F(\bar{Y}_{\hat{t}})=\vec{D}F(\bar{Y}_{\hat{t}}) \left(\bar{Y}_t-\bar{Y}_{\hat{t}}\right) +R_F(t,\hat{t}).
\end{equation}
\begin{lemma}\label{lem:RemainderF}
Let Assumption \ref{ass:1} hold. Then, for $t-\hat{t}<1$ and $ p\geq 1 $
$$
\Lpnorm{R_{F}(t,\hat{t})}   <C_p  (t-\hat{t})
$$
\end{lemma}

\begin{proof}
By defining a variable $Z:= \bar{Y}_{\hat{t}}+\theta \left(\bar{Y}_t-\bar{Y}_{\hat{t}}\right)$, we have the remainder term in  integral form,   
$$R_F(t,\hat{t})=\int_0^1{(1-\theta)D^2 F\left(Z\right) \Big( \bar{Y}_t-\bar{Y}_{\hat{t}},\bar{Y}_t-\bar{Y}_{\hat{t}}\Big) d \theta  }.
$$
Taking the norm, using Assumption \ref{ass:1}, Jensen's and H\"older's inequality
\begin{align}
 \Lpnorm{R_F(t,\hat{t})} 
   &\leq  \int_0^1{(1-\theta)\LpnormR{D^2 F\left(Z \right) \Big( \bar{Y}_t-\bar{Y}_{\hat{t}},\bar{Y}_t-\bar{Y}_{\hat{t}}\Big)} d \theta  }  \nonumber \\
    &\leq  \int_0^1\LpnormR{ \norm{  D^2 F\left(Z \right)}_{T_3}  \norm{\bar{Y}_t-\bar{Y}_{\hat{t}}}^2 } d \theta  \nonumber  \\
        &  \leq K \LpnormR{ \Big(  1+\norm{\bar{Y}_{\hat{t}}}^q +\norm{\bar{Y}_t}^q   \Big)    \norm{\bar{Y}_t-\bar{Y}_{\hat{t}}}^2 }  \nonumber \\
 &  \leq K \LtwopnormR{   1+\norm{\bar{Y}_{\hat{t}}}^q +\norm{\bar{Y}_t}^q   }   \LtwopnormR{\norm{\bar{Y}_t-\bar{Y}_{\hat{t}} }^2 }  \nonumber \\
          &  = K \LtwopnormR{   1+\norm{\bar{Y}_{\hat{t}}}^q +\norm{\bar{Y}_t}^q   }   \Lfourpnorm{\bar{Y}_t-\bar{Y}_{\hat{t}}  }^2 \nonumber \\
          & \leq C_p (t-\hat{t}) ,
\end{align}
where for the final inequality  Lemma \ref{lem:orderOfIncrementsInLp} is applied.
\end{proof}

Thirdly we define a remainder term $R$ from an \Ito-Taylor expansion of $F(\bar{Y}_t)$ near $F(\bar{Y}_{\hat{t}})$
\begin{equation} \label{eq:RemainderR}
R(t,\hat{t})=  F( \bar{Y}_t)  -   F( \bar{Y}_{\hat{t}})  -\vec{D}F(\bar{Y}_{\hat{t}}) \left( \sum _{i=1}^m   B_i \bar{Y}_{\hat{t}} \right) \left(  W_t^i-W^i_{\hat{t}}\right).
\end{equation}
The following lemma determines the order of the term $R$.
\begin{lemma}\label{lem:RemainderR}
Let Assumption \ref{ass:1} hold. Then for $ p\geq 2 $,
$$
\Lpnorm{R(t,\hat{t})}<C_p (t-\hat{t}).
$$
\end{lemma}

\begin{proof}
The substitution of \eqref{eq:timeDifference} into \eqref{eq:deterministicTaylor} gives us
 $$
F(\bar{Y}_t)  - F( \bar{Y}_{\hat{t}})=\vec{D} F(\bar{Y}_{\hat{t}}) \left( \SG{t,\hat{t}}-I\right)\bar{Y}_{\hat{t}}  + (t-\hat{t})\vec{ D}  F( \bar{Y}_{\hat{t}})  \SG{t,\hat{t}} \tilde{F}(\bar{Y}_{\hat{t}})+R_F(t,\hat{t}).
 $$
 Therefore we can re-express $R$ as follows
 $$
 R(t,\hat{t})=\vec{D} F(\bar{Y}_{\hat{t}})\left(  R_{GBM}(t,\hat{t})+ (t-\hat{t})     \tilde{F}(\bar{Y}_{\hat{t}}) \right)+ R_F(t,\hat{t}) 
 $$
Similar to \cite[Lemma 3.6]{wang}, we can determine the order of each term above.

Taking the norm and using Lemmas \ref{lem:RemainderGBM} and \ref{lem:RemainderF} with Assumption \ref{ass:1} we find
\begin{align}
 \Lpnorm{R(t,\hat{t})} & \leq \Lpnorm{ \vec{D} F(\bar{Y}_{\hat{t}} )R_{GBM}(t,\hat{t}) } + \vert t-\hat{t}\vert  \Lpnorm{ \vec{D} F(\bar{Y}_{\hat{t}}) \tilde{F}(\bar{Y}_{\hat{t}})  }\nonumber  \\ &+\Lpnorm{ R_F(t,\hat{t})}  \leq C_p (t-\hat{t}).
\end{align}
\end{proof}

\subsection{Strong convergence}\label{sec:main}
We now state and prove our main convergence result.
\begin{theorem}
\label{thm:convergence}
Let Assumption \ref{ass:1} hold.
Let $X_t$  be the solution of \eqref{eq:EqAFB} and $\bar{Y}_{t}$ the solution from \expint{TamedEI0}, \eqref{eq:tamedEI0}. Then 
there  exists a family of real numbers $C_{p,T} \leq 0$, independent of $\Delta t$ with $p \geq 1$  such that 
$$
\left( \eval{\sup_{t \in [0,T]} \norm{X_t-\bar{Y}_t}^p} \right)^{\frac{1}{p}} \leq  C_{p,T} \Delta t.
$$
\end{theorem}  
 
 \begin{proof}
We prove the result for $p\geq 4$ and then
from H\"{o}lder's inequality, we can conclude  the desired result   for $1 \leq p <4$.

 By considering \eqref{eq:cts} the difference between the exact and numerical solution is given by
\begin{equation}
 d  (X_t- \bar{Y}_t) = \left( A ( X_t-\bar{Y}_t)+ F( X_t)-\SG{t,\hat{t}}  \tilde{F}( \bar{Y}_{\hat{t}})  \right) dt +  \sum _{i=1}
  ^m    B_i (X_t-\bar{Y}_t ) dW_t^i.
\end{equation}
Applying the \Ito formula we get
 \begin{align*}
  d  ( \norm{X_t- \bar{Y}_t} ^2 )=&\left(  2  \innerproduct{X_t- \bar{Y}_t}{ A ( X_t-\bar{Y}_t)}+2 \innerproduct{X_t- \bar{Y}_t}{F( X_t)-\SG{t,\hat{t}}  \tilde{F}( \bar{Y}_{\hat{t}})  }    \right. \\ 
&   \left.  + \sum _{i=1} ^m  \norm{ B_i (X_t-\bar{Y}_{t}) }^2 \right)dt+
  2\left( \sum _{i=1}^m \innerproduct{X_t- \bar{Y}_t}{  B_i (X_t-\bar{Y}_{t}) } \right) dW_t^i.
 \end{align*}
 
 By the Cauchy--Schwarz inequality, and using that the error at $t=0$ is zero, we find
  for all $t \in [0,T]$
 \begin{align} \label{eq:bound1}
 \norm{X_t- \bar{Y}_t} ^2 & \leq   2 \norm{A}  \int_0^t \norm{X_s- \bar{Y}_{s}} ^2 ds   + 2 \int_0 ^t  \innerproduct{X_s- \bar{Y}_{s}}{F( X_s)-\SG{s,\hat{s}}  \tilde{F}( \bar{Y}_{\hat{s}})  }ds \nonumber \\
  &+\sum _{i=1}^m  \norm{ B_i} ^2 \int_0^t \norm{ X_s-\bar{Y}_{s}  }^2 ds \nonumber  +2 \sum _{i=1}^m  \int_0^t  \innerproduct{X_s- \bar{Y}_{s}}{  B_i (X_s-\bar{Y}_{s} ) }  dW_s^i \nonumber  \\
 & := \int_0^tI_1(s,\hat{s}) ds+2 \int_0^t I_2 (s,\hat{s})ds+\int_0^tI_3(s,\hat{s}) ds + 2 \sum _{i=1}^m \int_0^t I_4^{(i)}(s,\hat{s})  dW_s^i.
 \end{align}

 For the term  $I_2$, we add and subtract the terms $F( \bar{Y}_{s})$  and $\tilde{F}( \bar{Y}_{\hat{s}})$  to get
\begin{align}
 I_2 (s,\hat{s}) = & \innerproduct{X_s- \bar{Y}_{s}}{F( X_s)-\SG{s,\hat{s}}  \tilde{F}( \bar{Y}_{\hat{s}}  )  } \nonumber \\
=&\innerproduct{X_s- \bar{Y}_{s}}{F( X_s)-F( \bar{Y}_{s})} 
+ \innerproduct{X_s- \bar{Y}_{s}}{F( \bar{Y}_{s})- \tilde{F}( \bar{Y}_{\hat{s}})} \nonumber \\
& + \innerproduct{X_s- \bar{Y}_{s}}{\tilde{F}( \bar{Y}_{\hat{s}})-\SG{s,\hat{s}}  \tilde{F}(  \bar{Y}_{\hat{s}})  } \nonumber\\
=: & I_{2,1}(s,\hat{s})+I_{2,2}(s,\hat{s})+I_{2,3}(s,\hat{s}).
\end{align}
The one sided global Lipschitz of $F$ gives 
\begin{equation}
I_{2,1}(s,\hat{s}) \leq  K \norm{X_s- \bar{Y}_{s}}^2. 
\end{equation}
Noting that
$$
\tilde{F}( \bar{Y}_{\hat{s}})=F( \bar{Y}_{\hat{s}})-\Delta t  \frac{F( \bar{Y}_{\hat{s}}) \norm{F( \bar{Y}_{\hat{s}})}}{1+\Delta t \norm{F( \bar{Y}_{\hat{s}})}}
$$
yields
\begin{align}
I_{2,2} (s,\hat{s})& =\innerproduct{X_s- \bar{Y}_{s}}{F( \bar{Y}_{s})- F( \bar{Y}_{\hat{s}})} 
+   \innerproduct{X_s- \bar{Y}_{s}}{ \Delta t   \frac{F( \bar{Y}_{\hat{s}}) \norm{F( \bar{Y}_{\hat{s}})}}{1+\Delta t \norm{F( \bar{Y}_{\hat{s}})}}  } \nonumber \\
&  \leq \innerproduct{X_s- \bar{Y}_{s}}{F( \bar{Y}_{s})- F( \bar{Y}_{\hat{s}}}+\dfrac{1}{2} \norm{X_s- \bar{Y}_{s}}^2 +  \dfrac{1}{2} \Delta t ^2   \norm{F( \bar{Y}_{\hat{s}})}^4
\end{align}
where the Cauchy--Schwarz inequality and the inequality $2ab \leq a^2 + b^2$  for  the second term are applied. 
Returning to \eqref{eq:bound1}, we  have
  \begin{align} \label{eq:bound2}
 \norm{X_t- \bar{Y}_t} ^2 \leq &\left(2 \norm{A}+ \sum _{i=1}^m   \norm{ B_i} ^2 +1+2K  \right) \int_0^t \norm{X_s- \bar{Y}_{s}} ^2 ds  \nonumber \\ 
 &+ 2\int_0 ^t  \innerproduct{X_s- \bar{Y}_{s}}{    F(\bar{Y}_{s})-F(\bar{Y}_{\hat{s}})  }  ds + 2 \int_0 ^t  \innerproduct{X_s- \bar{Y}_{s}}{ \tilde{F}( \bar{Y}_{\hat{s}})-\SG{s,\hat{s}}  \tilde{F}( \bar{Y}_{\hat{s}})  }  ds   \nonumber  \\
 & + \Delta t ^2  \int_0 ^t      \norm{F(\bar{Y}_{\hat{s}}}^4 ds
   +2  \sum _{i=1}^m  \int_0^t  \innerproduct{X_s- \bar{Y}_{s}}{  B_i (X_s-\bar{Y}_{s} ) }    dW_s^i  .
 \end{align}
 
By taking supremum, we have 
  \begin{align} \label{eq:bound3}
 \sup_{s \in [0,T]} \norm{X_s- \bar{Y}_{s}} ^2 \leq & \left(2 \norm{A}+ \sum _{i=1}^m   \norm{ B_i} ^2 +1+2K  \right)  \int_0^T \norm{X_s- \bar{Y}_{s}} ^2 ds \nonumber \\ 
& +  2 \sup_{s \in [0,T]}  \int_0 ^s   \innerproduct{X_r-\bar{Y}_r}{F(\bar{Y}_{r})    -  F( \bar{Y}_{\hat{r}})  }  dr
 \nonumber \\
  &+2 \sup_{s \in [0,T]}  \int_0 ^s   \innerproduct{X_r-\bar{Y}_r}{ \tilde{F}( \bar{Y}_{\hat{r}})-\SG{r,\hat{r}}  \tilde{F}( \bar{Y}_{\hat{r}})   }   dr
 + \Delta t ^2  \int_0 ^T      \norm{F( \bar{Y}_{\hat{s}}}^4 ds \nonumber \\
 &  +  2 \sup_{s \in [0,T]}  \sum _{i=1}^m  \int_0^s   \innerproduct{X_r- \bar{Y}_r}{  B_i (X_r-\bar{Y}_r ) }   dW_r^i  .
 \end{align} 
Working in the space $L^{p/2} (\Omega, \mathbb{R})$  for $p \geq 4$, using that $\LphalfnormR{\norm{X}^2}=\Lpnorm{X}^2$ for a random variable in $\mathbb{R}^d$, we have
 
 \begin{align} \label{eq:bound4}
 & \LphalfnormR{\sup_{s \in [0,T]} \norm{X_s- \bar{Y}_{s}} ^2}  \nonumber \\
 &\leq \left(2 \norm{A}+ \sum _{i=1}^m   \norm{ B_i} ^2 +1+2K  \right) \int_0^T \Lpnorm{X_s- \bar{Y}_{s}} ^2 ds  \nonumber \\ 
 &+ 2  \LphalfnormR{  \sup_{s \in [0,T]}  \int_0 ^s  \innerproduct{X_r- \bar{Y}_r}{F( \bar{Y}_r)    -  F( \bar{Y}_{\hat{r}})  }ds } \nonumber \\
  &+ 2  \LphalfnormR{  \sup_{s \in [0,T]}  \int_0 ^s  \innerproduct{X_r- \bar{Y}_r}{  \tilde{F}( \bar{Y}_{\hat{r}})-\SG{r,\hat{r}}  \tilde{F}( \bar{Y}_{\hat{r}} (\hat{r})) }ds } \nonumber \\
 &+ \Delta t ^2  \int_0 ^T      \Ltwopnorm{F( \bar{Y}_{\hat{s}}}^4 ds 
   +  2 \LphalfnormR{\sup_{s \in [0,T]} \sum _{i=1}^m  \int_0^s  \innerproduct{X_r- \bar{Y}_r}{  B_i (X_r-\bar{Y}_r ) }  dW_r^i}.
 \end{align} 
 We consider the last term separately. By  applying the triangle inequality, the inequality in  \cite[Lemma 3.7]{hutzenthaler2012},   Cauchy--Schwarz inequality, H\"{o}lder inequality  and arithmetic-geometric mean inequality, we have
 \begin{align*}
 2 & \LphalfnormR{\sup_{s \in [0,T]}  \sum _{i=1}^m  \int_0^s  \innerproduct{X_r- \bar{Y}_r}{  B_i (X_r-\bar{Y}_r ) }  dW_r^i}\\
& \leq 2  \sum _{i=1}^m   \LphalfnormR{\sup_{s \in [0,T]} \int_0^s  \innerproduct{X_r- \bar{Y}_r}{  B_i (X_r-\bar{Y}_r ) }  dW_r^i}\\
   & \leq p  \sum _{i=1}^m   \left(   \int_0^T \LphalfnormR{ \innerproduct{X_s- \bar{Y}_{s}}{  B_i (X_s-\bar{Y}_{s} ) }}^2  ds \right) ^{1/2} \\
          & \leq  p  \sum _{i=1}^m    \left(  \int_0^T \LphalfnormR{ \norm{X_s- \bar{Y}_{s}}\norm{{  B_i (X_s-\bar{Y}_{s} ) }}}^2  ds \right) ^{1/2} \\
          & \leq  p \sum _{i=1}^m   \left(   \int_0^T \Lpnorm{X_s- \bar{Y}_{s}}^2 \Lpnorm{  B_i (X_s-\bar{Y}_{s} ) }^2  ds \right) ^{1/2} \\
          & \leq   \sup_{s \in [0,T]}  \Lpnorm{X_s- \bar{Y}_{s}}  p \sum _{i=1}^m  \left(    \int_0^T \Lpnorm{  B_i (X_s-\bar{Y}_{s} ) }^2  ds \right) ^{1/2} \\ 
          & \leq \frac{1}{4}   \LpnormR{  \sup_{s \in [0,T]} \norm{X_s- \bar{Y}_{s}}} ^2+  p^2 m \sum _{i=1}^m  \int_0^T \Lpnorm{  B_i (X_s-\bar{Y}_{s} ) }^2  ds.
\end{align*}
  Substitution of this result into \eqref{eq:bound4} results in

 \begin{align} \label{eq:bound5}
 & \frac{3}{4}   \LphalfnormR{\sup_{s \in [0,T]} \norm{X_s- \bar{Y}_{s}} ^2}  = \frac{3}{4} \LpnormR{\sup_{s \in [0,T]} \norm{X_s- \bar{Y}_{s}} }^2  \nonumber \\
  & \leq \left(2 \norm{A}+ \sum _{i=1}^m   \norm{ B_i} ^2(1+p^2m) +1+2K  \right) \int_0^T \Lpnorm{X_s- \bar{Y}_{s}} ^2 ds  \nonumber \\ 
& + \Delta t ^2  \int_0 ^T      \Ltwopnorm{F(\bar{Y}_{\hat{s}}}^4 ds    
 + 2  \LphalfnormR{  \sup_{s \in [0,T]}  \int_0 ^s  \innerproduct{X_r- \bar{Y}_r}{F(\bar{Y}_r)    -  F( \bar{Y}_{\hat{r}})  }dr } \nonumber \\
& +2  \LphalfnormR{  \sup_{s \in [0,T]}  \int_0 ^s  \innerproduct{X_r- \bar{Y}_r}{  \tilde{F}( \bar{Y}_{\hat{r}})-\SG{r,\hat{r}}  \tilde{F}( \bar{Y}_{\hat{r}})  }dr }   \nonumber  \\
  &:= T_1 + T_2 + T_3 + T_4.
\end{align} 

We need to expand the third and fourth terms. Starting with third term by applying Cauchy--Schwarz inequality and arithmetic-geometric mean inequality, we have
\begin{multline}
T_3=2  \LphalfnormR{  \sup_{s \in [0,T]}  \int_0 ^s  \innerproduct{X_r- \bar{Y}_r}{F(\bar{Y}_r) -F( \bar{Y}_{\hat{r}})}dr }\\
\leq T_{3,1} + 2  \LphalfnormR{  \sup_{s \in [0,T]}  \int_0 ^s  \innerproduct{X_r- \bar{Y}_r}{  R(r,\hat{r}) }dr }\\
\leq T_{3,1}+ \int_0^T \Lpnorm{X_r- \bar{Y}_{r}} ^2 dr +\int_0^T \Lpnorm{R(r,\hat{r})}^2 dr
\end{multline}

where  $R$ is the remainder term explicitly defined in \eqref{eq:RemainderR} and 
$$
T_{3,1}=2  \LphalfnormR{  \sup_{s \in [0,T]}  \int_0 ^s \innerproduct{X_r- \bar{Y}_r}{     \vec{D}F(\bar{Y}_{\hat{r}}) \left( \sum _{i=1}^m   B_i   \bar{Y}_{\hat{r}} \right) \left(  W_r^i-W^i_{\hat{r}}\right) }dr }\ .
$$ 

To determine the order of $T_{3,1}$
we expand $X_r$ and $\bar{Y}_{r}$ as follows
$$
X_r=\SG{r,\hat{r}} \left(X_{\hat{r}}+\int_{\hat{r}}  ^r \SG{\tau,\hat{r} } ^{-1} F(X_{\tau}) d \tau \right),\qquad
\bar{Y}_{r}=\SG{r,\hat{r}} \left( \bar{Y}_{\hat{r}}+\int_{\hat{r}}  ^r  \tilde{F}( \bar{Y}_{\hat{r}}) d \tau \right)
$$
and take the difference to obtain 
\begin{multline}
X_r- \bar{Y}_{r}=\SG{r,\hat{r}} \left(X_{\hat{r}}-\bar{Y}_{\hat{r}}\right) +  \SG{r,\hat{r}} \int_{\hat{r}}  ^r \left[   \SG{\tau,\hat{r} } ^{-1} F(X_{\tau})- F(X_{\tau})\right]  d\tau\\ +   \SG{r,\hat{r}} \int_{\hat{r}}  ^r    \left[ F(X_{\tau})- \tilde F(\bar{Y}_{\hat{r}})\right]  d\tau.
\end{multline}
By  adding and subtracting $ \SG{r,\hat{r}} F ( X_{\hat{r}} )$, we have 
\begin{multline}\label{eq:1st_term_inner_product}
X_r- \bar{Y}_r=\SG{r,\hat{r}} \int_{\hat{r}}  ^r    \left[ F(X_{\tau})-  F( X _{\hat{r}})\right]  d\tau+\SG{r,\hat{r}} \int_{\hat{r}}  ^r \left[   \SG{\tau,\hat{r} } ^{-1} F(X_{\tau})- F(X_{\tau})\right]  d\tau\\+\SG{r,\hat{r}} \left(X_{r}-\bar{Y}_{\hat{r}}+(r-\hat{r})F( X_{\hat{r}})- (r-\hat{r})\tilde {F}(  \bar{Y}_{\hat{r}})  \right).
\end{multline}

Expanding $T_{3,1}$ using \eqref{eq:1st_term_inner_product}  we have
\begin{align*}
T_{3,1} \leq& 2\LphalfnormR{  \sup_{s \in [0,T]}  \int_0 ^s  \innerproduct{\SG{r,\hat{r}} \int_{\hat{r}}  ^r    \left[ F(X_{\tau})- F(X_{\hat{r}})\right]  d\tau}{      \Theta _{r,\hat{r}}}dr } \\
 &+ 2\LphalfnormR{  \sup_{s \in [0,T]}  \int_0 ^s  \innerproduct{\SG{r,\hat{r}} \int_{\hat{r}}  ^r \left[   \SG{\tau,\hat{r} } ^{-1} F(X_{\tau})-  F(X_{\tau} )\right]  d\tau}{    \Theta _{r,\hat{r}}  }dr}\\
 &+ 2\LphalfnormR{  \sup_{s \in [0,T]}  \int_0 ^s  \innerproduct{\SG{r,\hat{r}} \zeta _{r,\hat{r}}}{      \Theta _{r,\hat{r}}}dr}  \\
& + 2\LphalfnormR{  \sup_{s \in [0,T]}  \int_0 ^s  \innerproduct{\SG{r,\hat{r}} \left( X_{\hat{r}}-\bar{Y}_{\hat{r}} \right) }{ \Theta _{r,\hat{r}}} dr}\\
:=& T_{3,1,1}+T_{3,1,2}+T_{3,1,3}+T_{3,1,4}
\end{align*}
where 
$$
\Theta _{r,\hat{r}}=\vec{D}F( \bar{Y}_{\hat{r}}) \left( \sum _{i=1}^m   B_i   \bar{Y}_{\hat{r}} \right) \left(  W_r^i-W^i_{\hat{r}}\right)
$$
and 
$$
  \zeta_{r,\hat{r}}= (r-\hat{r})F(X_{\hat{r}})- (r-\hat{r})\tilde{F}(  \bar{Y}_{\hat{r}}).
 $$
    
 By definition and  by order property of Brownian increment in $L^p(\Omega,\mathbb{R}^d)$ and polynomial growth condition on $\vec{D}F$ and bounded moments of numerical solutions, it is straightforward to write 
 \begin{equation} \label{eq:ThetaOrder}
\Lpnorm{\vec{\Theta} _{r,\hat{r}}} \leq C_p \sqrt{r-\hat{r}}.
\end{equation}
By applying Cauchy--Schwarz and H\"older inequalities  respectively, we have  
 $$
  T_{3,1,1} \leq     \int_0 ^T  \int_{\hat{r}} ^r  \LtwopnormR{ \norm{\SG{r,\hat{r}}}}    \Ltwopnorm{ F(X_{\tau})- F(X_{\hat{r}} }  \Lpnorm{    \Theta _{r,\hat{r}} } d\tau   dr  \leq C_{p,T} \Dt^2
 $$
 where we have used Lemma \ref{lem:SGBoundedness}-ii, Lemma \ref{lem:orderOfIncrementsInLp}, equation \eqref{eq:ThetaOrder} and  boundedness  of   the   terms stated in Lemmas \ref{lem:boundedFexact} and  \ref{lem:boundedFnumeric}.
 In a similar way, but considering the inequality  (ii)  given in Lemma \ref{lem:GBM_F_order},  we have  
 $$
 T_{3,1,2} \leq 2   
 \int_0 ^T  \Lpnorm{\SG{r,\hat{r}} \int_{\hat{r}}  ^r \left[   \SG{\tau,\hat{r} } ^{-1} F(X_{\tau})- F(X_{\tau})\right]  }\Lpnorm{    \Theta _{r,\hat{r}}  } d\tau dr\leq C_{p,T} \Delta t ^2.
 $$
 
By adding and subtracting, $\zeta _{r,\hat{r}}$ to the first term of the inner product in $T_{3,1,3}$, we  have 
\begin{align*}
T_{3,1,3} &=  2\LphalfnormR{  \sup_{s \in [0,T]}  \int_0 ^s  \innerproduct{\SG{r,\hat{r}} \zeta _{r,\hat{r}}}{      \Theta _{r,\hat{r}}}dr}  \\
          &\leq  2\LphalfnormR{  \sup_{s \in [0,T]}  \int_0 ^s  \innerproduct{\SG{r,\hat{r}} \zeta _{r,\hat{r}}- \zeta _{r,\hat{r}}}{      \Theta _{r,\hat{r}}}dr} +2\LphalfnormR{  \sup_{s \in [0,T]}  \int_0 ^s  \innerproduct{ \zeta _{r,\hat{r}}}{      \Theta _{r,\hat{r}}}dr} \\
 & \leq    2 \int_0 ^T  \Lpnorm{ \SG{r,\hat{r}} \zeta _{r,\hat{r}}- \zeta _{r,\hat{r}}  } \Lpnorm{    \Theta _{r,\hat{r}}  } dr + 2  \LphalfnormR{  \sup_{s \in [0,T]}  \int_0 ^s  \innerproduct{ \zeta _{r,\hat{r}}}{      \Theta _{r,\hat{r}}}dr}        
\end{align*}
The first term above is of order $\Dt^2$ as in $T_{3,1,1}$ and $T_{3,1,2}$. The second term, which is denoted by $J_4$ in \cite{wang},  is determined as  $O(\Dt^2)$ in  \cite[Eq. (3.45) to Eq. (3.56)]{wang} by use of Doob's maximal inequality and a Burkholder-Davis-Gundy type inequality for discrete-time martingales. Therefore, we find $T_{3,1,1}< C \Dt ^2$. 
 
We now examine $T_{3,1,4}$. By adding and subtracting $X_{\hat{r}}-\bar{Y}_{\hat{r}}$   to the first term of the inner product in $T_{3,1,4}$, we  have 
\begin{align*} 
T_{3,1,4} & \leq 
  2\LphalfnormR{  \sup_{s \in [0,T]}  \int_0 ^s  \innerproduct{ \SG{r,\hat{r}} \left( X_{\hat{r}}-\bar{Y}_{\hat{r}} \right)-\left( X_{\hat{r}}-\bar{Y}_{\hat{r}} \right) }{ \Theta _{r,\hat{r}}} dr} \\
& \qquad + 2\LphalfnormR{  \sup_{s \in [0,T]}  \int_0 ^s  \innerproduct{  \left( X_{\hat{r}}-\bar{Y}_{\hat{r}} \right) }{ \Theta _{r,\hat{r}}} dr} \\
&=:  T_{3,1,4,1}+ T_{3,1,4,2}.
\end{align*}
Applying H\"{o}lder's inequality and the arithmetic-geometric mean inequality to $T_{3,1,4,1}$ 
\begin{align*}
T_{3,1,4,1} & \leq \int_0 ^T  \Lpnorm{  X_{\hat{r}}-\bar{Y}_{\hat{r}  } }     \Lpnorm{ \left(\SG{r,\hat{r}} -I \right)^\intercal  \Theta _{r,\hat{r}}}    dr \\
& \leq    \sup_{s \in [0,T]}   \Lpnorm{  X_s-\bar{Y}_s }  \int_0 ^T  \Lpnorm{ \left(\SG{r,\hat{r}} -I \right)^\intercal  \Theta _{r,\hat{r}}}    dr  \\
& \leq  \frac{1}{4}  \sup_{s \in [0,T]}   \Lpnorm{  X_s-\bar{Y}_s } ^2 +C \Dt^2.
\end{align*}
The term $T_{3,1,4,2}$, which is denoted by $J_5$ in \cite[Lemma 3.7]{wang}, is also bounded in a similar way to $T_{3,1,4,1}$. 

Finally putting these bounds together for the third term in \eqref{eq:bound5}, we have
\begin{equation}\label{eq:T3Order}
T_3 \leq  C_{p,T} \Dt^2+  \frac{1}{4} \sup_{s \in [0,T]}  \Lpnorm{X_s- \bar{Y}_{s} }^2+ C_{p,T} \int_0^T \Lpnorm{X_s- \bar{Y}_{s}} ^2 ds.
\end{equation}
The term $T_4$ in \eqref{eq:bound5} remains to be considered. By the expansion \eqref{eq:R_GBM} 
, we see that 
\begin{multline*}
T_4=2  \LphalfnormR{  \sup_{s \in [0,T]}  \int_0 ^s  \innerproduct{X_r- \bar{Y}_r}{  \tilde{F}( \bar{Y}_{\hat{r}})-\SG{r,\hat{r}}  \tilde{F}( \bar{Y}_{\hat{r}})  }dr }\\
\leq 
2  \LphalfnormR{  \sup_{s \in [0,T]}  \int_0 ^s \innerproduct{X_r- \bar{Y}_r}{  -\left( \sum _{i=1}^m   B_i   \bar{Y}_{\hat{r}} \right) \left(  W_r^i-W^i_{\hat{r}}\right) }dr }\\
+ 2  \LphalfnormR{  \sup_{s \in [0,T]}  \int_0 ^s  \innerproduct{X_r- \bar{Y}_r}{ - R_{GBM}(r,\hat{r}) }dr }.
\end{multline*}
Denoting the first term by $T_{4,1}$ and using the Cauchy-Schwarz inequality
\begin{multline}
T_4
\leq T_{4,1}+ \int_0^T \Lpnorm{X_s- \bar{Y}_{s}} ^2 ds +\int_0^T \Lpnorm{R_{GBM}(s,\hat{s})}^2ds
\end{multline}
where  $R_{GBM}$ is remainder term which is explicitly defined in \eqref{eq:remainderGBM} and its order is determined in Lemma \ref{lem:RemainderGBM}.
If the  expansion of $X_r-\bar{Y}_r$ given in  \eqref{eq:1st_term_inner_product} is substituted in  $T_{4,1}$, by following same steps for $T_3$, we obtain exactly the same bound as for $T_3$.
Inserting these bounds for $T_3$ and $T_4$      into \eqref{eq:bound5}, and applying the Gronwall inequality completes the proof for $p \geq 4$. 
\end{proof}
\section{Extensions of the scheme and numerical results}
\subsection{An Adaptive GBM based scheme with  TamedEI0 as a backstop method}
\label{sec:backstop}
   The drawback of taming  methods is  the use  of  the modified  drift functions at every step,
   even in the case that  the numerical solution remains small. 
   As noted in the introduction, in the recent literature, this issue is handled  in  several ways.
   We examine the 
   adaptive methods proposed by Kelly et al. \cite{kelly2016adaptive,KellyMilsteinArxiv} that employ both usual and tamed time stepping schemes. If the numerical solution stays in the region determined by the admissibility condition, usual time stepping scheme is applied. Otherwise, a tamed version of the time stepping scheme is chosen (termed a backstop method).
  
   In our numerical experiments in Section \ref{Sec:Num}, we observe an adaptive strategy  employing EI0 introduced in \cite{utkuLord}  and the \expint{TamedEI0}  \eqref{eq:tamedEI0}  as a backstop performs  remarkably well.

The scheme we examine numerically is given by 
\begin{equation}  \label{eq:adaptiveEI0}
  Y_{n+1}^{adp}=\SG{t_{n+1},t_n}\left( Y_n^{adp}+ F(Y_n^{adp}) h_{n+1}\mathbbm{1}_{\lbrace  h_{n+1} >h_{\min}\rbrace}  + 
  \tilde{F}(Y_n^{adp})\mathbbm{1}_{\lbrace h_{n+1}= h_{\min} \rbrace }   \right) 
\end{equation}
where
$$
 h_{n+1} (Y_n)= \max \left\lbrace  h_{\min}, \min \left\lbrace  h_{\max},\dfrac{ h_{\max}}{\norm{Y_n} }\right\rbrace \right\rbrace 
$$ and  $t_n=\sum_{i=1}^{n} h_i$  with fixed ratio  $\dfrac{h_{\max}}{h_{\min}}=\rho.$ 

To prove  convergence of the scheme \eqref{eq:adaptiveEI0}, one would need to follow the steps in \cite{kelly2016adaptive,KellyMilsteinArxiv}. 
The major issue for \eqref{eq:adaptiveEI0}  would be to establish suitable one-step error bounds as in \cite[Section 3]{KellyMilsteinArxiv} (rather than the final time bounds here) and this is the subject of future work.
  
\subsection{Application to a larger class of SDEs} \label{sec:LargerClass}

The appearance  of  non-linear $g_i:\mathbb{R}^d \to \mathbb{R}^d $ functions in diffusion coefficients in \eqref{eq:EqAFB} yields a  semi-linear SDE
\begin{equation} \label{eq:EqAFBg}
  dX_t= \left( A X_t+F(X_t) \right) dt + \sum _{i=1}
  ^m  \left(  B_i X_t+g_i(X_t)\right)  dW_t^i 
\end{equation}
where $X_{0}=\xi \in  \mathbb{R} ^d$ and $m$, $d \in
  \mathbb{N}$. We proposed the following scheme  in \cite{utkuLord} for globally Lipschitz $F$ and $g_i$ functions
  \begin{equation}  \label{eq:fullEI0}
    Y_{n+1}^N=\SG{t_{n+1},t_n}\left( Y_n+  \left( F
      (Y_n ^ N) - \sum _{i=1} ^m B_i g_i(Y_n ^N) \right)  \Delta t + \sum _{i=1} ^m g_i(\vec{Y}_n^N
      )\Delta W_{t_n} ^i \right)
  \end{equation} 
 and proved strong convergence of order $\frac{1}{2}$.

 When the drift term $F$ satisfies Assumption \ref{ass:1} it is natural to propose drift tamed versions of \eqref{eq:fullEI0} by utilizing the modified $\tilde{F}$ in \eqref{eq:taming} instead of $F$ in the above and would expect convergence rate of order $1/2$.

\subsection{Numerical Results}
\label{Sec:Num}
In this section, we illustrate convergence and compare the methods in this paper with recent methods appearing in the literature. In the following examples, the acronyms {\tt GBM} and {\tt AGBM} correspond to the the tamed GBM method \expint{TamedEI0} \eqref{eq:tamedEI0} with fixed step size and the adaptive version of GBM of \eqref{eq:adaptiveEI0} respectively introduced here.
First we examine some small dimensional SDEs and then the discretization of an SPDE.

For the SDE examples we denote the tamed Milstein of \cite{wang} by {\tt Mil} and the adaptive Milstein method of \cite{KellyMilsteinArxiv} by {\tt AMil}.
We let {\tt Pmil} denote the Projected Milstein of \cite{Beyn2017}. We take $2000$ realizations in each case and present error bars based on $20$ groups of $100$ on the root mean square error (RMSE). In each case we take $10^6$ steps to form a reference solution and (other than where an analytic solution was available) we used the Tamed Milstein method to construct this. For efficiency we plot the RMSE against an average clock time from the computations, termed cputime below.

In our final example, based on a spectral discretization of a stochastic PDE we compare the GBM based methods to a semi-implicit tamed and exponential based tamed methods and we examine the effectiveness and reduction in order of the scheme proposed in Section \ref{sec:LargerClass}.

\subsection{Ginzburg Landau Equation}
Consider the scalar SDE
\begin{equation}
\label{eq:GL}
dX_t= \left( - X_t +\frac{\sigma}{2}X_t - X_t^3 \right)  dt+\sqrt{\sigma} X_t dW(t)
\end{equation}
with exact solution as given in  \cite{kloeden2011} 
\begin{equation}
X_t=\frac{X_0e^{-t+\sqrt{\sigma} W(t)} }{\sqrt{1+2X_0^2 \int _0 ^t  e^{-2s+2\sqrt{\sigma} W(s)}ds}}.
\end{equation}
This equation is often used as a test equation, for example in \cite{hutzenthaler2011} to show the divergence of explicit Euler-Maruyama scheme. 
\begin{figure}[ht]
\begin{center}
    \includegraphics[width=0.495\textwidth,height=0.45\textwidth]{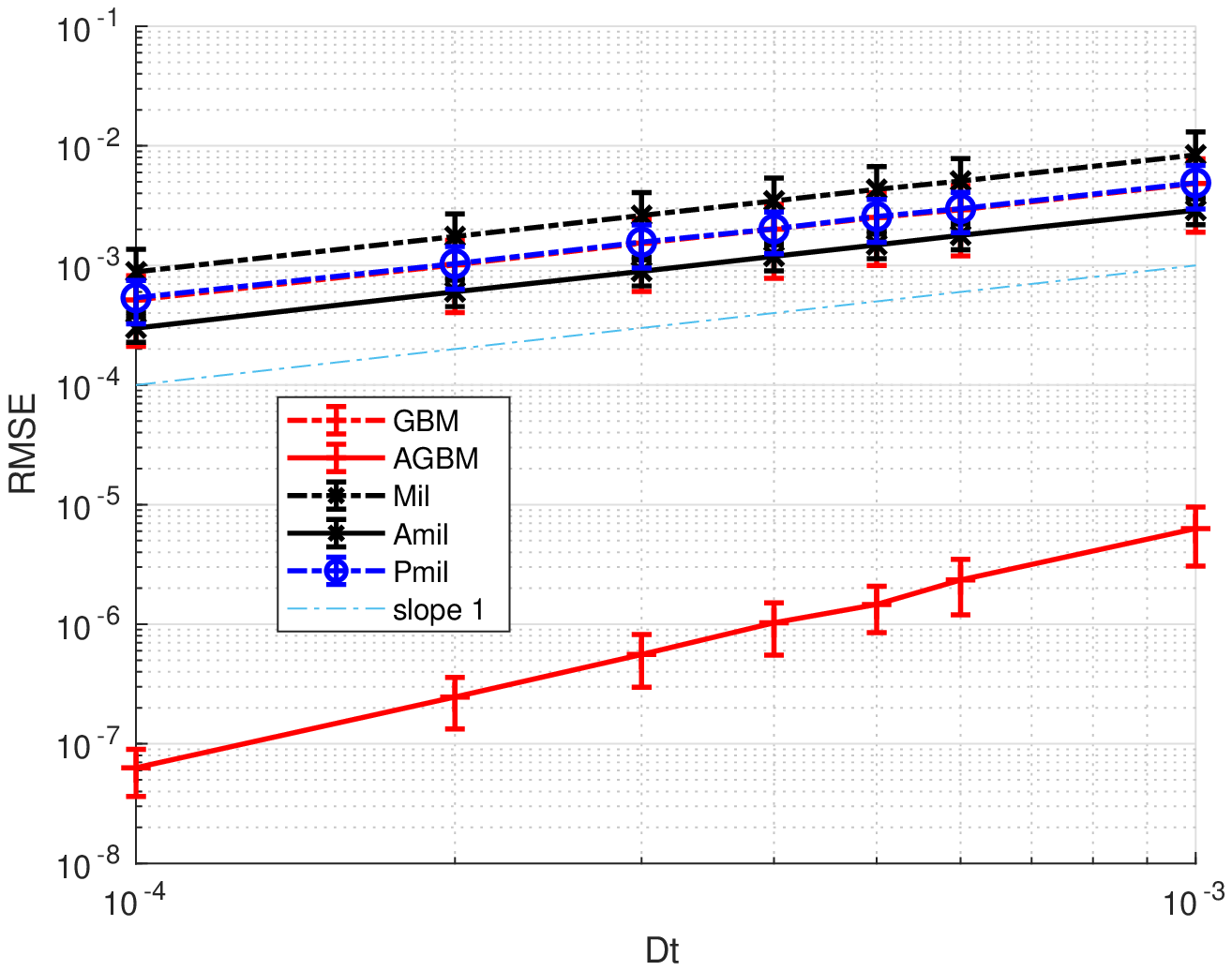}
    \includegraphics[width=0.495\textwidth,height=0.45\textwidth]{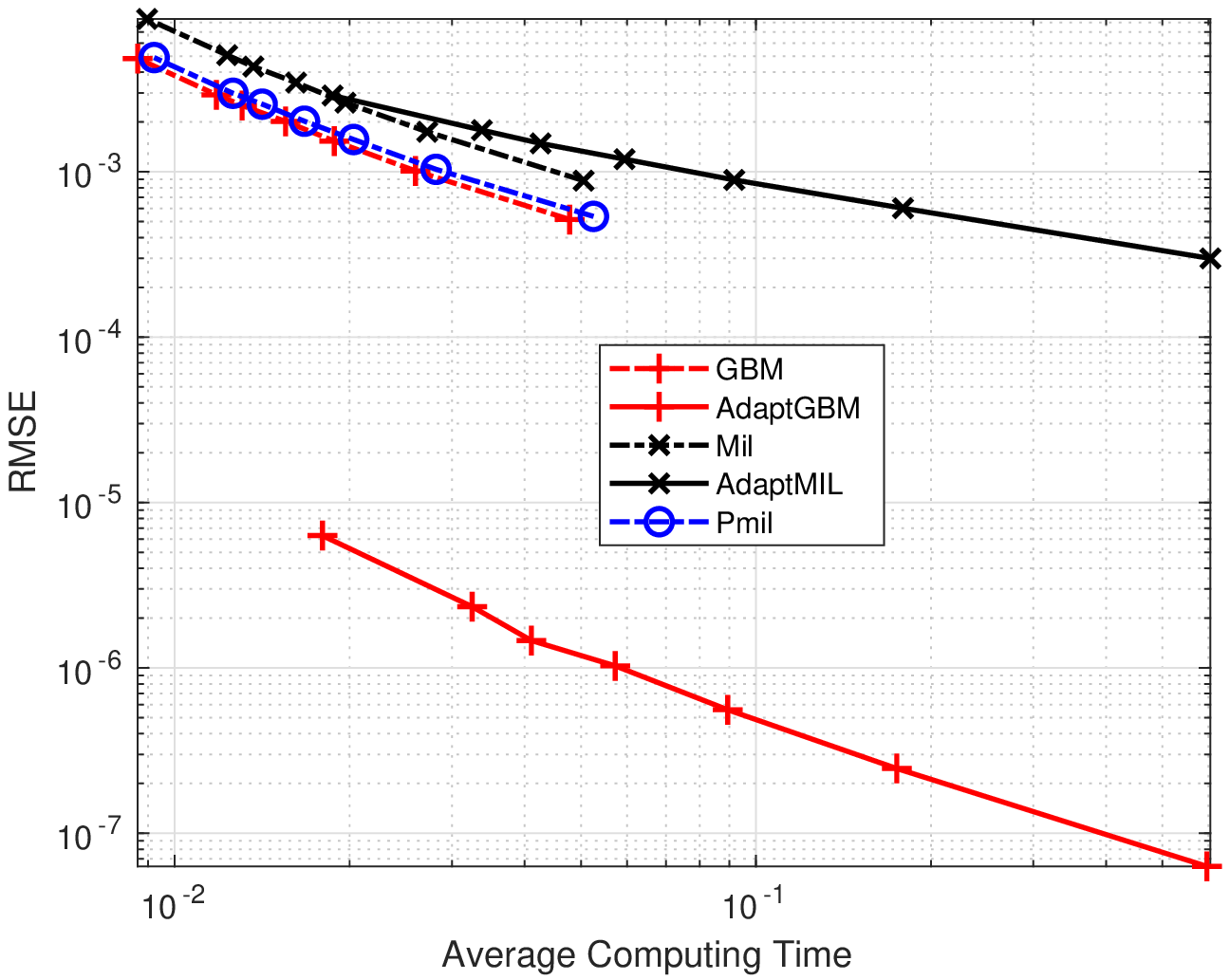}
\end{center}
\caption{ Stochastic Ginzburg-Landau Equation \eqref{eq:GL} with
  $\sigma=2$, $T=1$ (a) root mean square error (RMSE)
  against $\Dt$ and a reference line with slope 1.
  (b) RMSE against cputime.}  
\label{fig:Landau}
\end{figure}
In Fig. \ref{fig:Landau} (a) we see that {\tt GBM} is the most accurate of the fixed step methods and that {\tt AGBM} clearly has a much smaller error constant in this case. In (b) we see that {\tt AGBM} is clearly the most efficient overall followed by {\tt GBM}.

\subsection{HIV Internal Dynamics}
\begin{figure}[ht] \label{fig:HIV}
\begin{center}
\includegraphics[width=0.495\textwidth,height=0.45\textwidth]{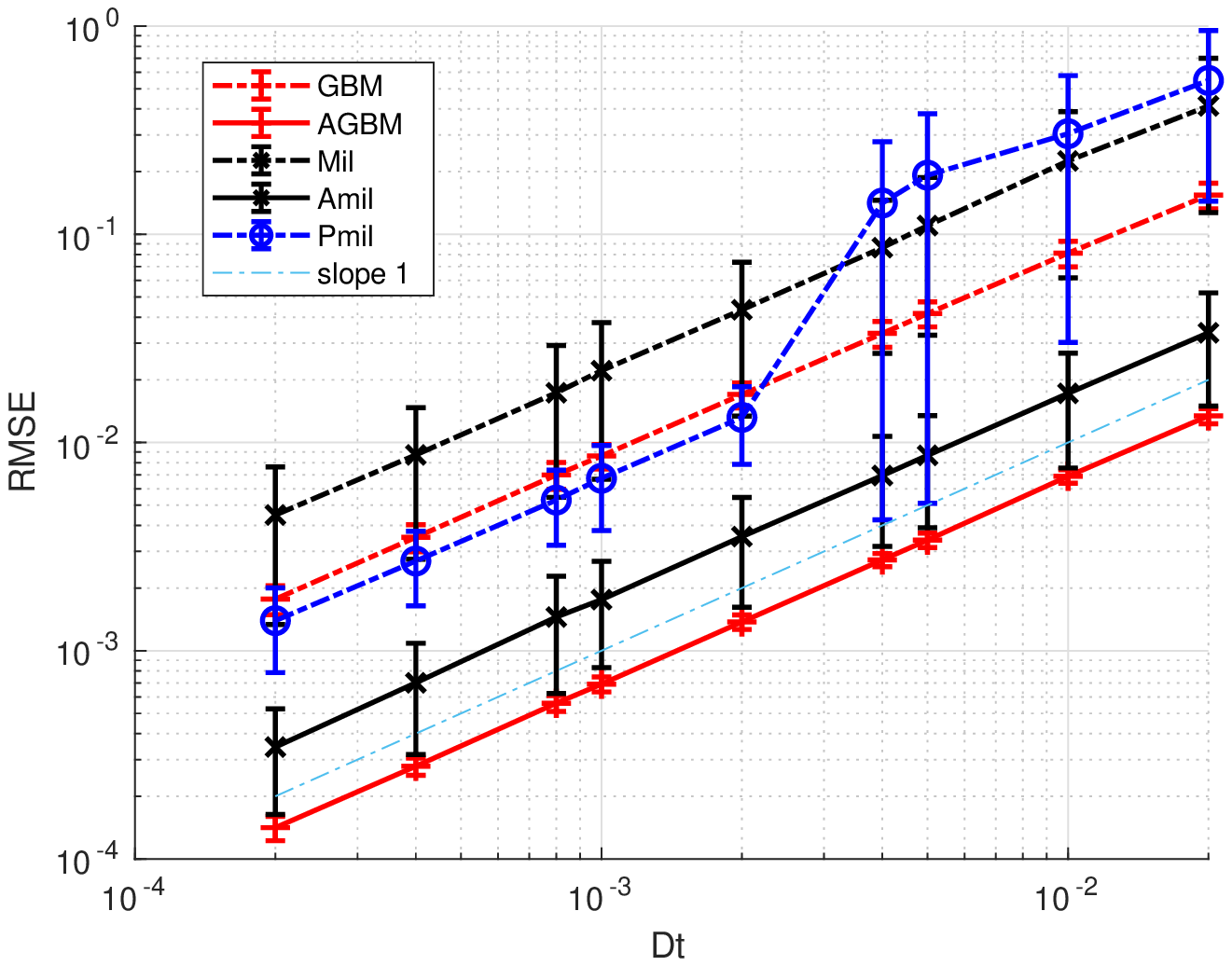}
\includegraphics[width=0.495\textwidth,,height=0.45\textwidth]{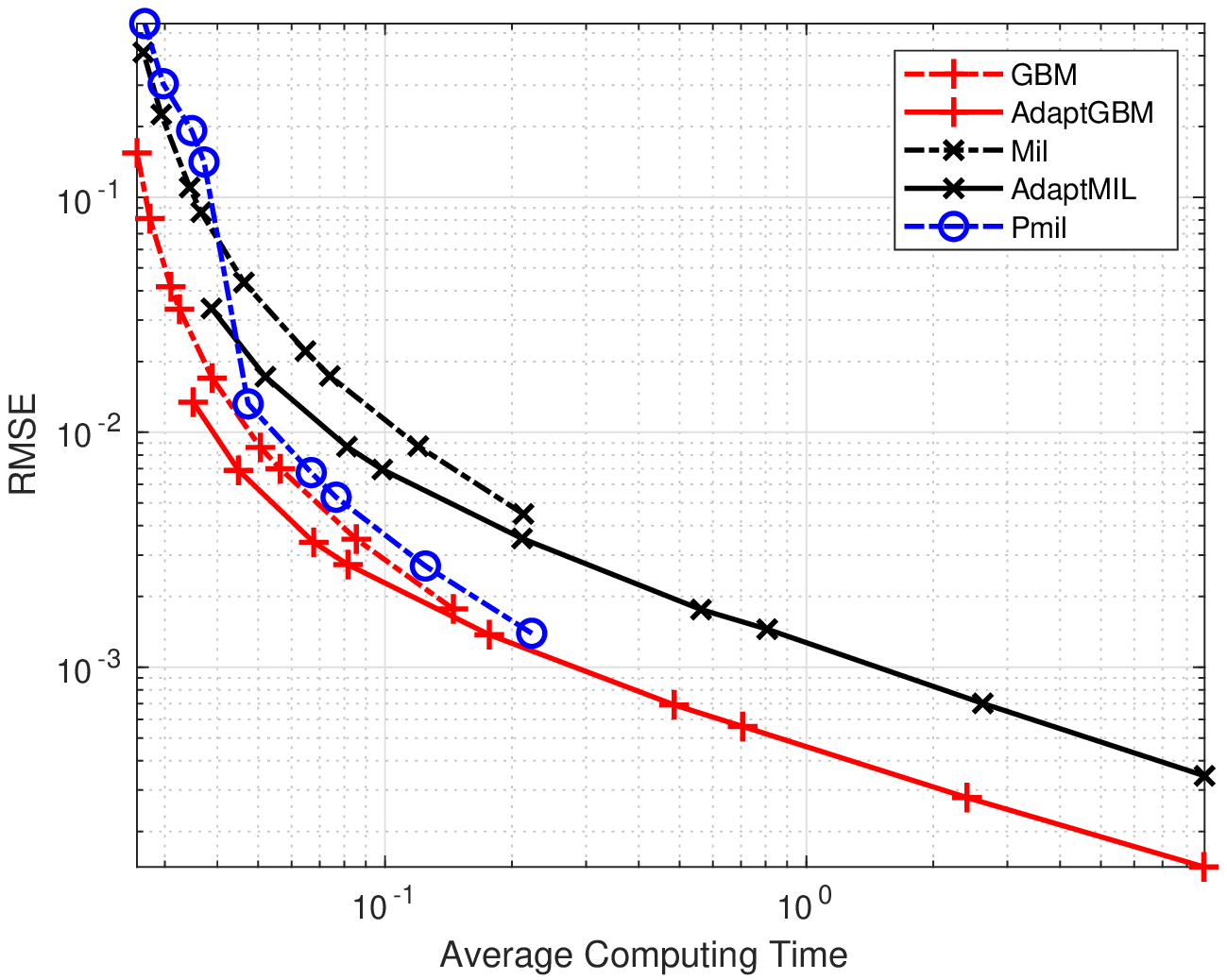}
\end{center}
\caption{Stochastic HIV Model \eqref{eq:HIV} with parameters $\lambda=3$, $\mu=2$, $k=0.5$, $\alpha=0.7$, $c=0.1$ and noise parameters  $\sigma_1=1.25$, $\sigma_2=0.09$, $\sigma_3=0.4$. In (a) we plot RMSE against $\Dt$ with a reference line with slope 1 and (b) RMSE against cputime.}
\end{figure}

Consider the SDE model for HIV internal dynamics introduced in \cite{sonner2016,MaoHiv2008}
\begin{align}\label{eq:HIV}
 dT&=(\lambda - \mu T -kTV)dt+\sigma_1TdW^1 \nonumber \\
 dI&=(kTV-\alpha I)dt+\sigma_2 I dW^1  \\ 
 dV&=(cI - \gamma V  -kTV)dt+\sigma_3 VdW^2  \nonumber     
\end{align}
where the parameters values taken are given in the caption to Fig. \ref{fig:HIV} and $T_0=0.5$, $I_0=0.7$ and $V_0=0.9$.
To construct the GBM based methods, we take $A=\text{diag}(-\mu,-\alpha,-\gamma)$  and 
$\sum_{i=1}^2 B_i=\text{diag}(\sigma_1,\sigma_2,\sigma_3)$. 
In Fig. \ref{fig:HIV} (a) we see the two adaptive methods are more accurate than the fixed step methods and that the overall {\tt AGBM} is the most accurate and from (b) the most efficient. {\tt GBM} performs well over a large range of step sizes and this leads to an efficient scheme in (b). We note that the projected method has a drop in the error as $\Dt$ is reduced as there is less need for the projection.
\subsection{Lotka--Volterra SDE}
\begin{figure}[ht]
\begin{center}
\includegraphics[width=0.495\textwidth,height=0.45\textwidth]{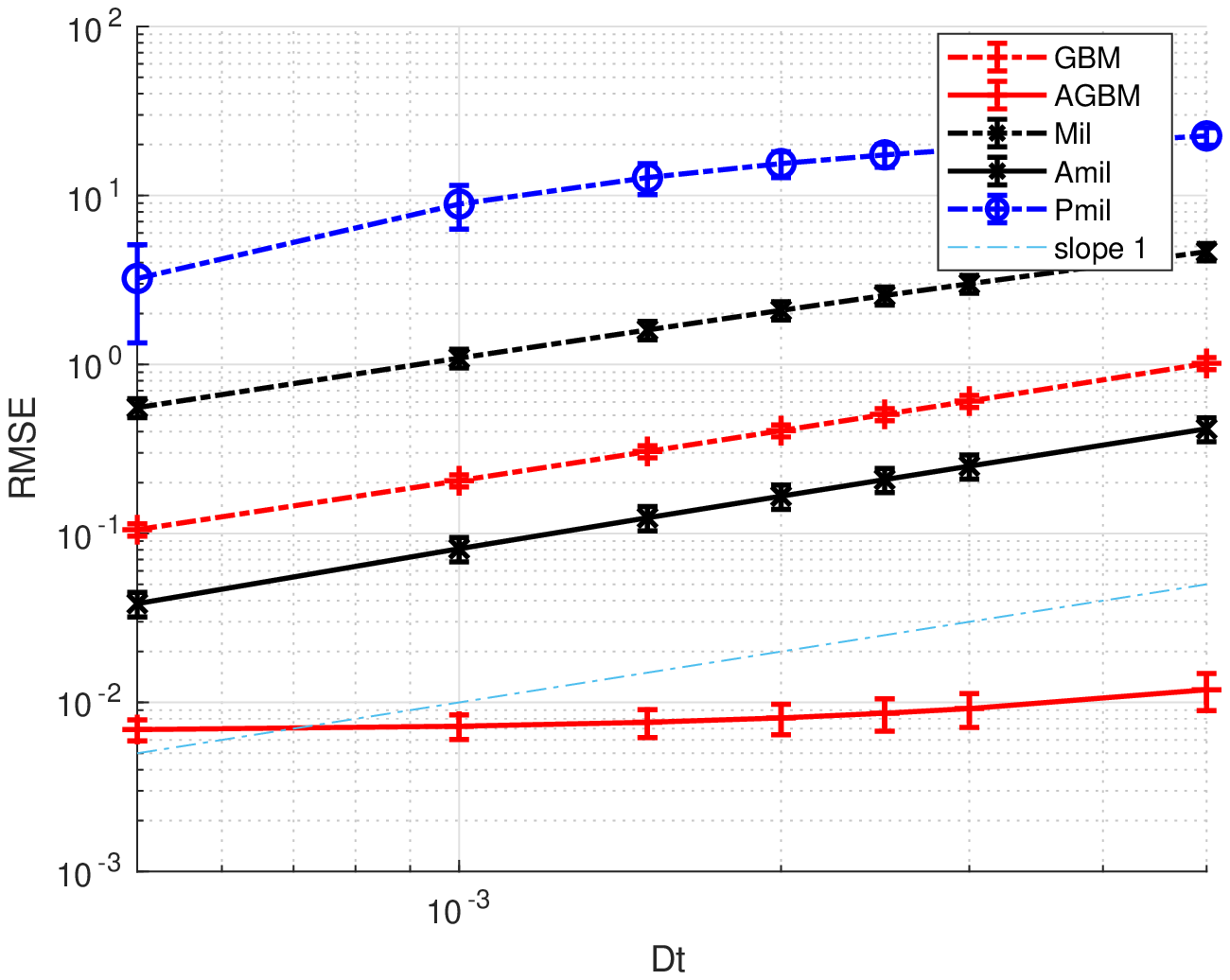}
\includegraphics[width=0.495\textwidth,height=0.45\textwidth]{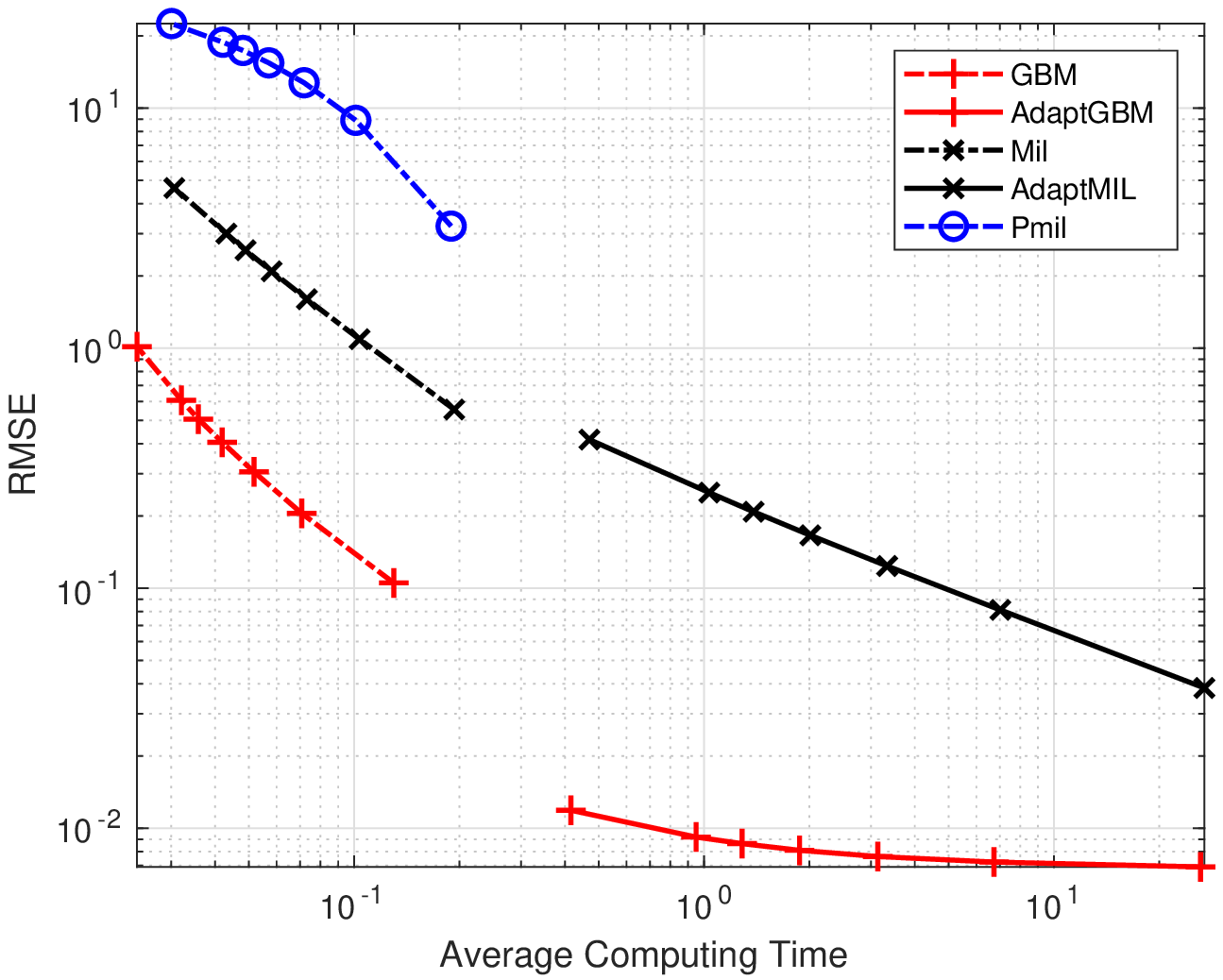}
\end{center}
\caption{Stochastic Lotka--Volterra  Model \ref{eq:LV} with parameters
 $\lambda=0.8 $,$\beta=0.15$,$\delta=0.75$,$\gamma=0.01$,$\sigma_1=\sigma_2=0.1$ with $T=1$. (a) RMSE against $\Dt$ and a reference line with slope 1 and 
  (b) RMSE  against cputime.}
   \label{fig:LV}
\end{figure}

We consider a stochastic 
Lotka--Volterra model as studied in  \cite{kelly2016adaptive}, 
\begin{align} \label{eq:LV}
dX&=X(\lambda - \beta Y)dt+ \sigma_1 X dW^1 \\
dY&=Y(\gamma X - \delta )dt+ \sigma_2 Y dW^2.
\end{align}
The parameter values taken are given in Fig. \ref{fig:LV} and for initial data we took $X_0=5$ and $Y_0=10$.
We constructed the GBM based methods by taking
$A=\text{diag}(\lambda,-\delta)$  and  
$\sum_{i=1}^2 B_i=\text{diag}(\sigma_1,\sigma_2)$.
In Fig. \ref{fig:LV} we see {\tt GBM}  has the smallest error amongst the fixed step methods and that {\tt AGBM} is the most accurate overall - indeed the tamed Milstein reference solution is not sufficiently accurate to observe further convergence. In this example the error in {\tt Pmil} is large and we are only just starting to see convergence as $\Dt$ is reduced. It is interesting to note that the fixed step {\tt GBM} in this case is more efficient than the adaptive {\tt Amil}.

\subsection{Model of tumor growth}
Consider the growth model for tumor cells under the influence of random perturbations  which is analysed in  \cite{tumorGrowth}  
\begin{equation} \label{eq:CM}
dp_t = \left( \lambda \ln{\frac{\mu}{p_t}}-G(v(t))\right)p_t dt
+\sigma p_t dW
\end{equation}
where $p_t$ the number of cancerous tumor cells and $v(t)=(1+\cos(t))^{-1}$ is the dose of the drug at  time t, $G(v(t))=k_1 v(t)(k_2+v(t))^{-1}$ is the destroying rate per tumor cell and time unit. 
To apply the GBM  based integrators, we take $A=0$,  $F(p)=\left( \lambda \ln{\frac{\mu}{p}}-G(v)\right)p$ and $B=\sigma$. The resulting scheme  becomes
$$
p_{n+1}=\exp(-\frac{1}{2}\sigma^2 \Dt+\sigma \Delta W_n   ) \left(p_n+  \Dt F(p_n)   \right).
$$
We take $p_0=0.8$ and $T=1$.
\begin{figure}
\begin{center}
\includegraphics[width=0.495\textwidth,height=0.45\textwidth]{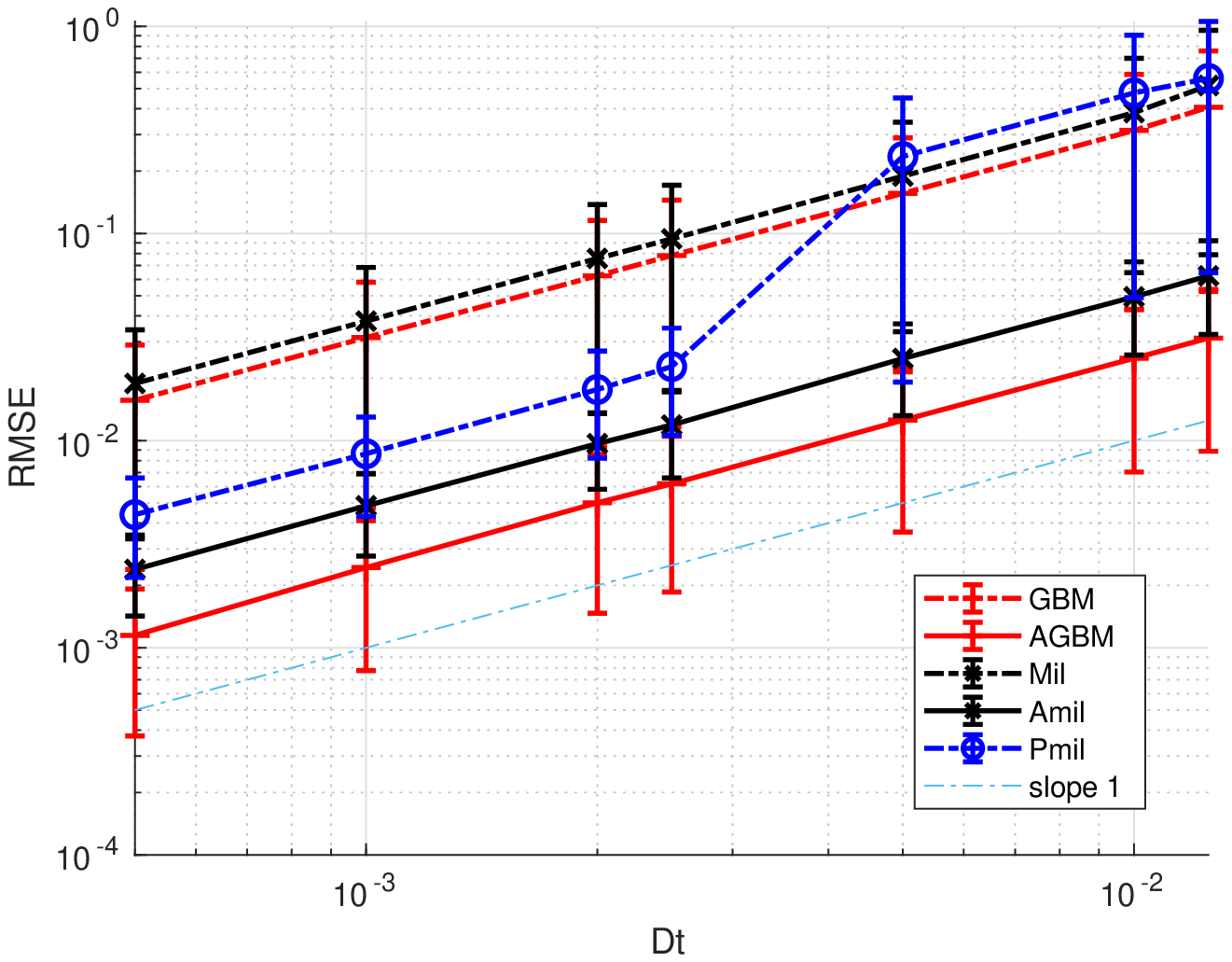}
\includegraphics[width=0.495\textwidth,height=0.45\textwidth]{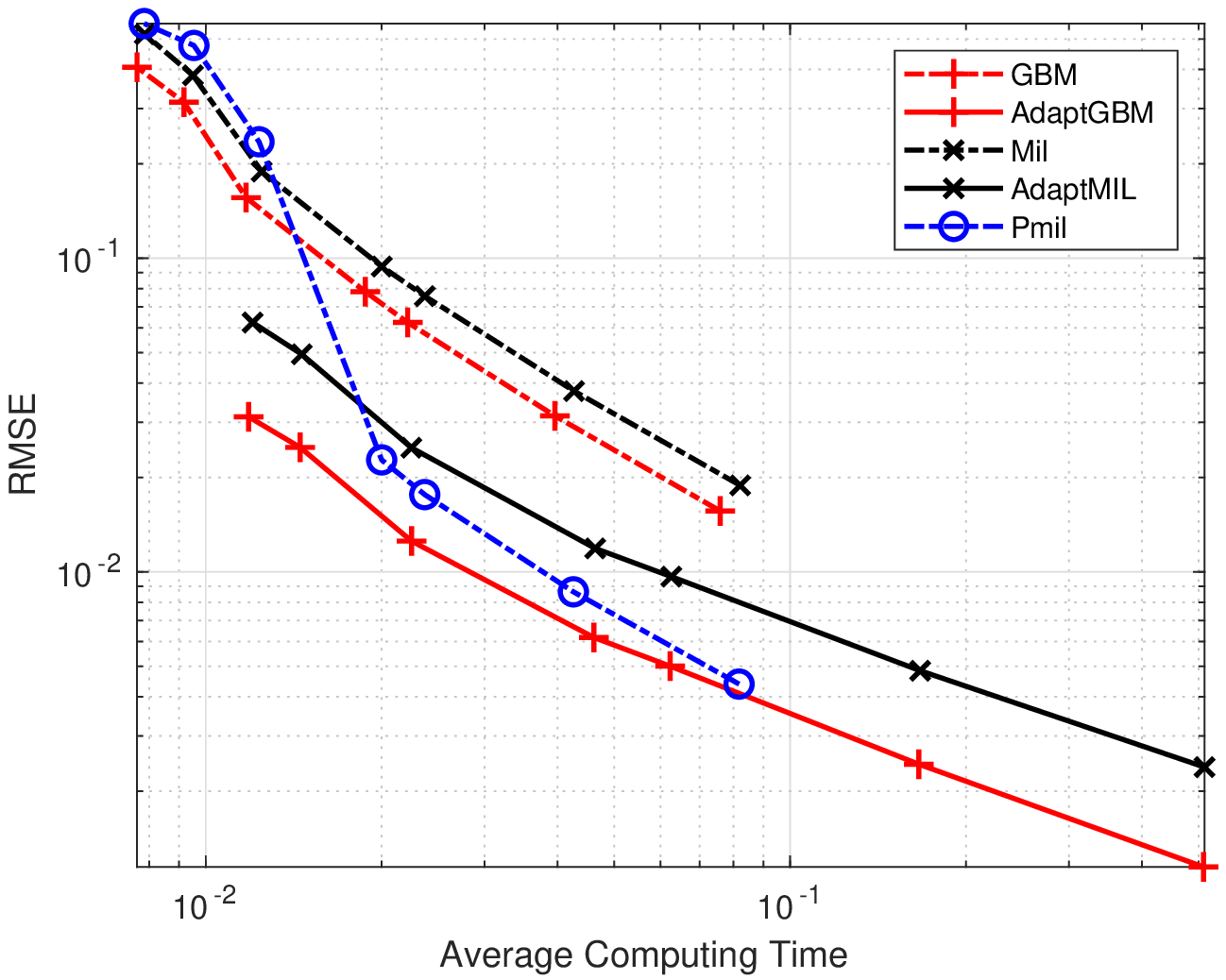}
\end{center}
\caption{Stochastic tumor growth model \ref{eq:CM} with parameters $\lambda=1$, $\mu=1$, $k_1=k_2=1$, $\sigma=1.5$.In (a) RMSE against $\Dt$ and reference line with slope 1 in (b) RMSE against cputime.}
\label{fig:CM}
\end{figure}  
 In Fig. \ref{fig:CM} we see that {\tt AGBM} is again the most accurate and efficient although in this example {\tt Pmil} is the leading fixed step method once in the asymptotic regime.

\subsection{Semi-Linear SPDE}
Consider the stochastic reaction-diffusion equation 
\begin{equation} \label{eq:SPDE}
du=\left[\varepsilon \frac{\partial ^2  u}{\partial x^2}+u-\gamma u^3
\right]dt+\left[ \beta u +\alpha \frac{1-u}{1+u^2} \right]dW,
\qquad u(x,0)=\sin(\pi x),
\end{equation}
with $x\in[0,1]$ subject to zero Dirichlet boundary conditions. We take $W$
to be a $Q$--Wiener process and let the covariance operator $Q$ have orthonormal eigenfunctions
$g_j(x)=\sqrt{2}sin(j\pi x)$ and eigenvalues $\nu_j =\frac{1}{j^2}$,
$j \in \mathbb{N}$, so that  
$$W = \sum_{j \in \mathbb{N}} \frac{1}{j} g_j \beta_j,$$
where $\beta_j$ are iid Brownian motions.

We applied the spectral Galerkin method in space, described for example in \cite{expmil}, with $d=128$ Fourier components to obtain a system of SDEs.
For this large system we took $M=500$ samples, taking groups of 20 to estimate the standard deviation in the root mean square error (RMSE) and $10^5$ steps to construct the reference solution.

When $\alpha=0$ this system of SDEs is of the form of
\eqref{eq:EqAFB} and we can apply \eqref{eq:tamedEI0} to the system of SDEs. For $\alpha\neq 0$ the SDEs are in the form of \eqref{eq:EqAFBg} and we can apply the tamed scheme in \ref{sec:LargerClass}.
For $\alpha=0$ 
we benchmark against a semi-implicit tamed Milstein {\tt Imp Mil}, a semi-implicit adaptive tamed Milstein {\tt Adaptive Imp Mil} and exponential versions {\tt Exp Mil} and {\tt Adaptive Exp Mil}. These schemes are chosen so that mean-square linear stability is not an issue and for SDEs they would all be first order.
For $\alpha \neq 0$ we compare to a semi-implicit EM method {\tt Imp EM} and an adaptive version {\tt Adaptive Imp EM} as well as exponential versions {\tt ETD} and {\tt Adaptive ETD}. Again mean-square linear stability is not an issue and for SDEs they would be order $1/2$.

\begin{figure}
\begin{center}
\includegraphics[width=0.495\textwidth,height=0.45\textwidth]{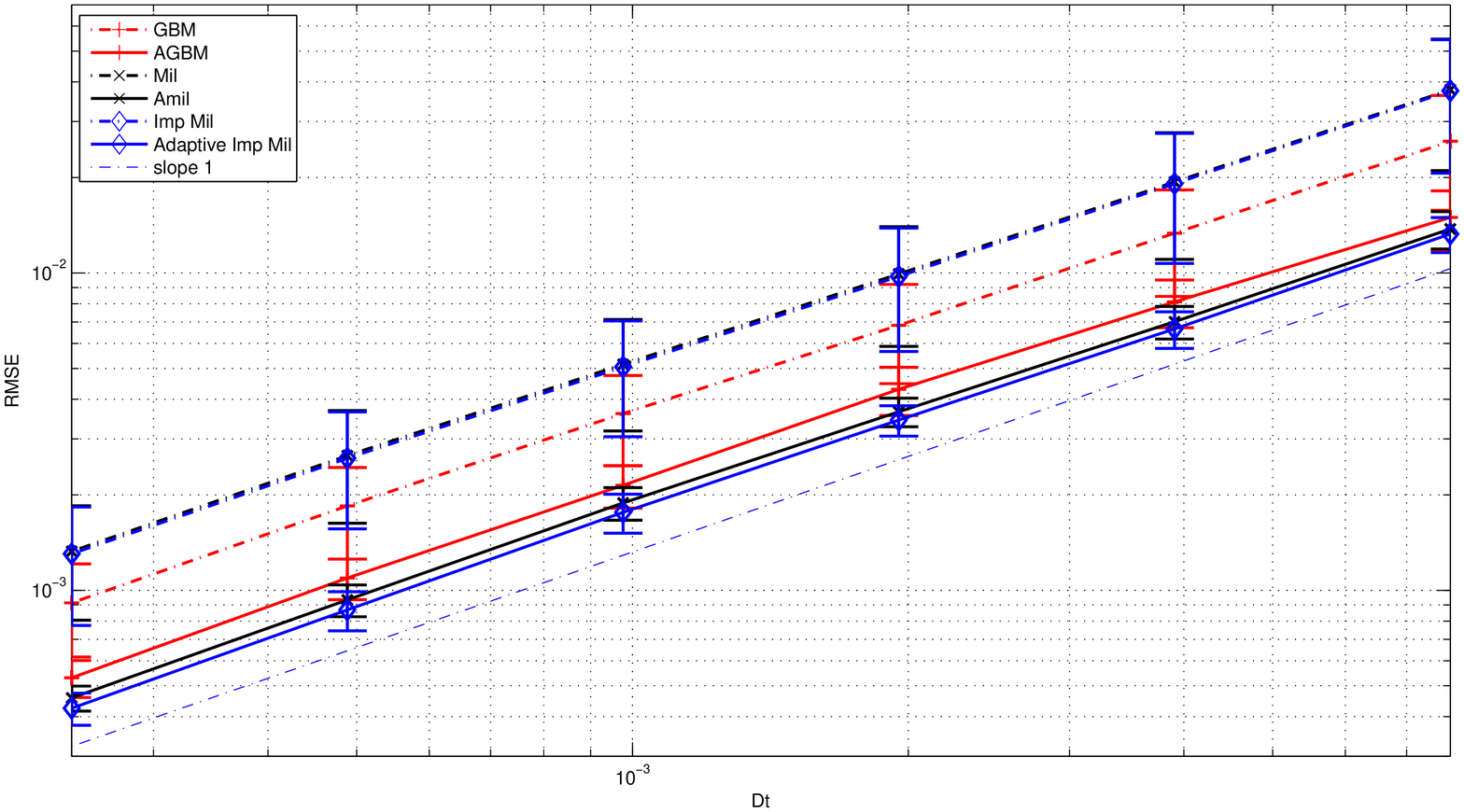}
\includegraphics[width=0.495\textwidth,height=0.45\textwidth]{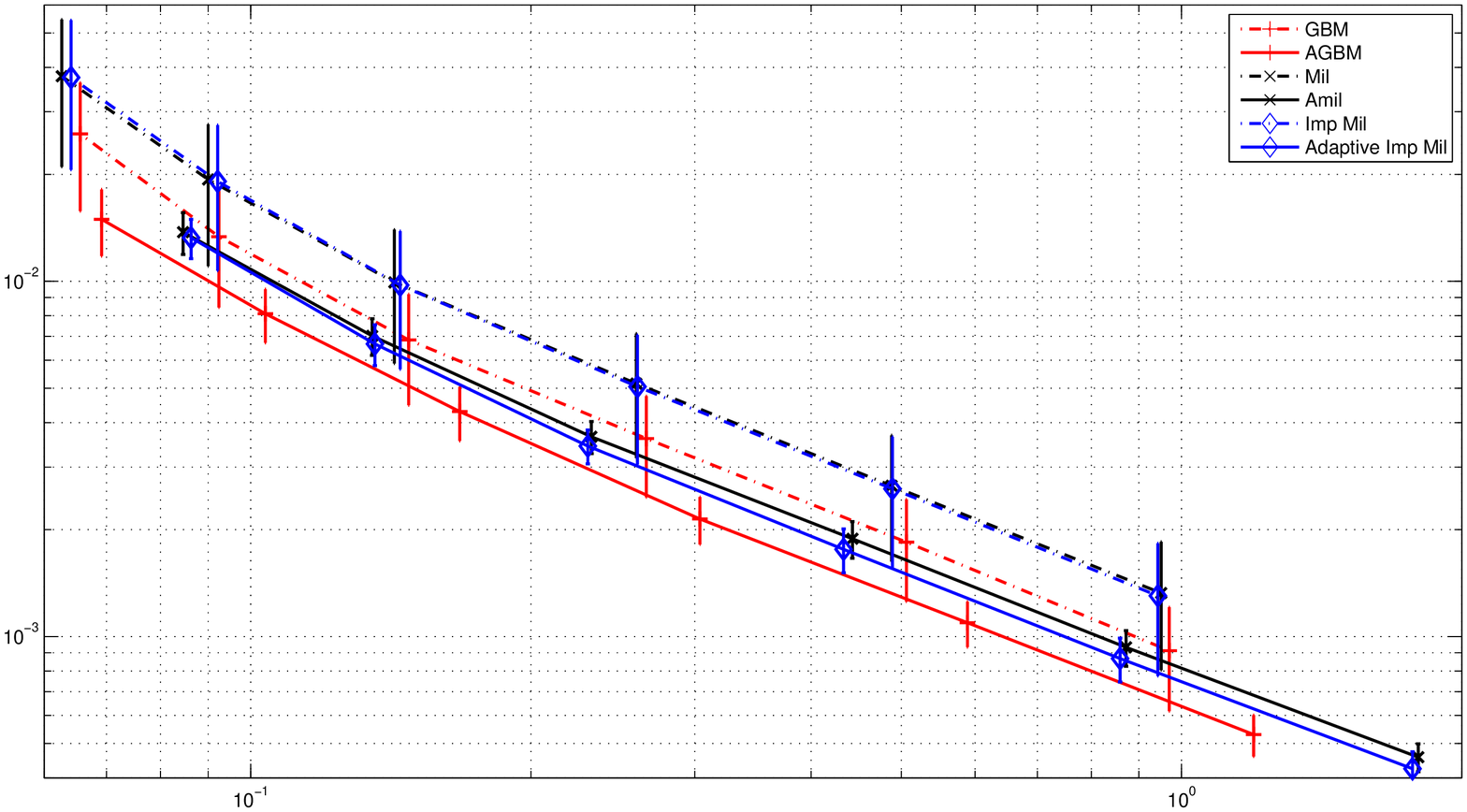}
\end{center}
\caption{SPDE \eqref{eq:SPDE} with  $\alpha=0$, $\beta=1$, $\gamma=1$ with $M=500$ samples (a) RMSE against $\Dt$ and reference line with slope 1 (b) RMSE against cputime.} 
\label{fig:alpha0beta1}
\end{figure}

In Fig. \ref{fig:alpha0beta1} we examine linear noise ($\alpha=0$)
and we observe convergence of order $1$ in all methods and see that {\tt Adaptive GBM} is the most efficient and {\tt GBM} the most efficient fixed step method.

\begin{figure}
\begin{center}
\includegraphics[width=0.495\textwidth,height=0.45\textwidth]{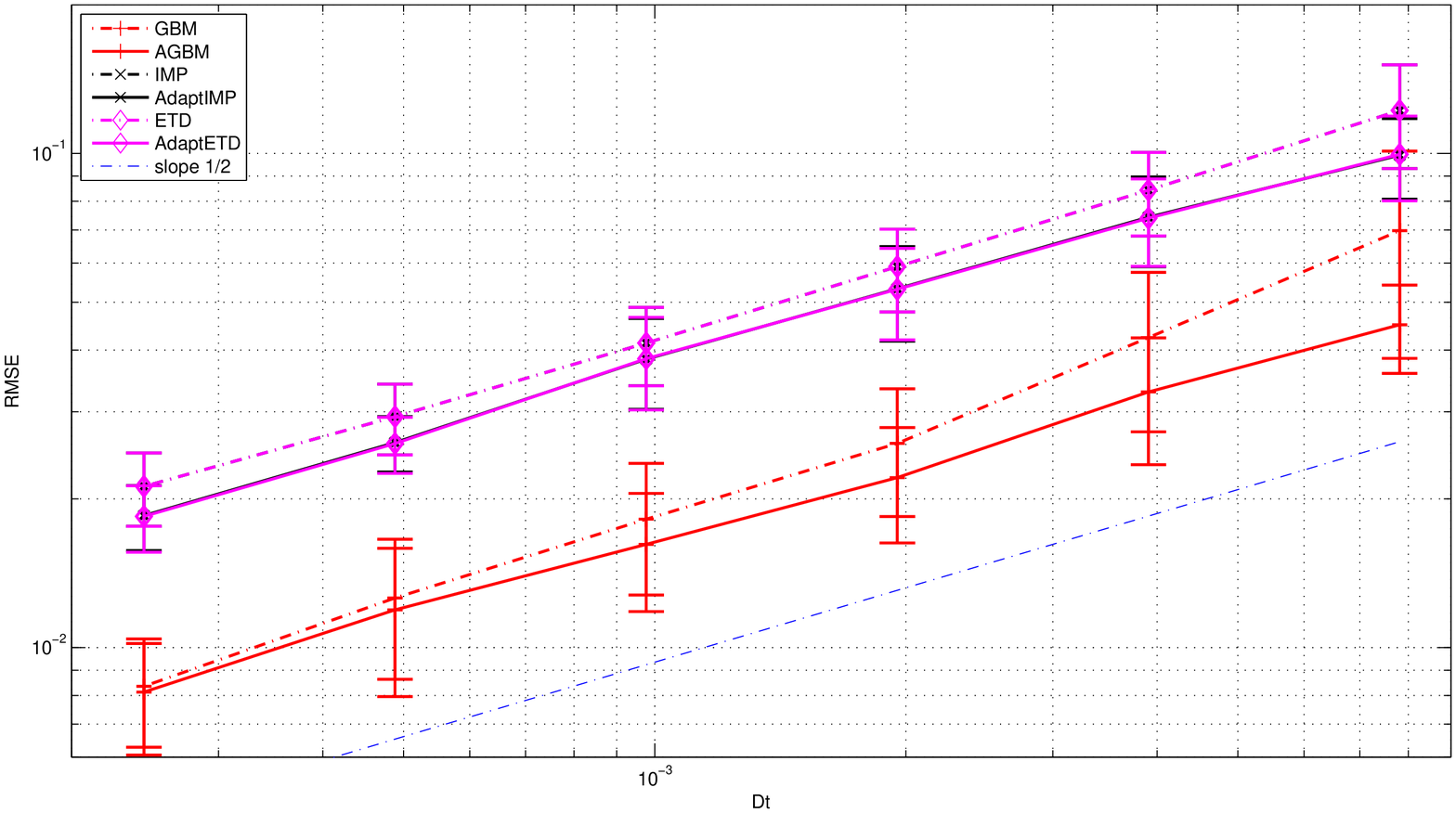}
\includegraphics[width=0.495\textwidth,height=0.45\textwidth]{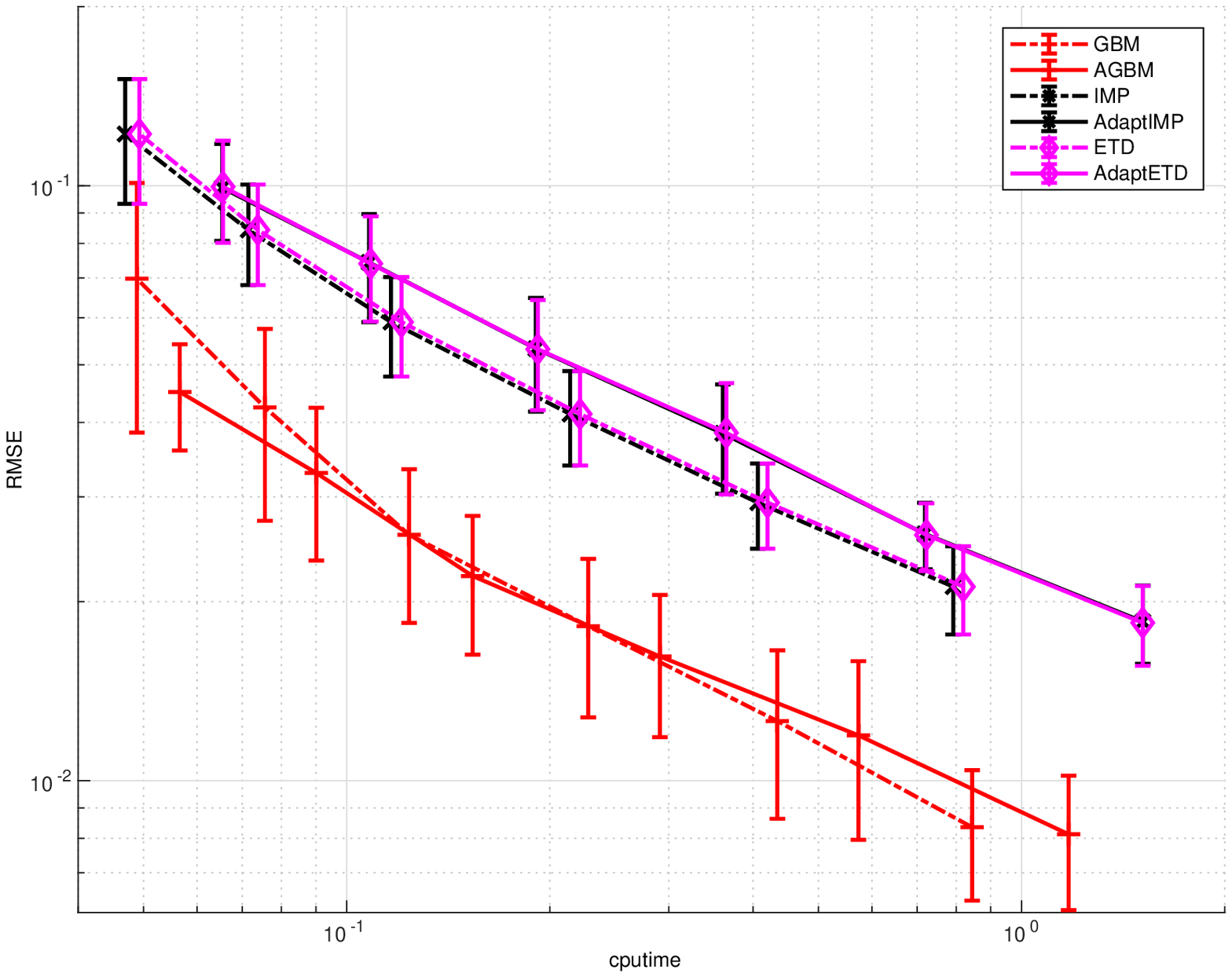}
\end{center}
\caption{SPDE \eqref{eq:SPDE} with  $\alpha=0.5$, $\beta=1$, $\gamma=0.25$ with $M=500$ samples (a) RMSE against $\Dt$ and a reference line with slope $\frac{1}{2}$ (b)  RMSE against cputime.} 
\label{fig:alpha05beta1}
\end{figure}

\begin{figure}
\begin{center}
\includegraphics[width=0.495\textwidth,height=0.45\textwidth]{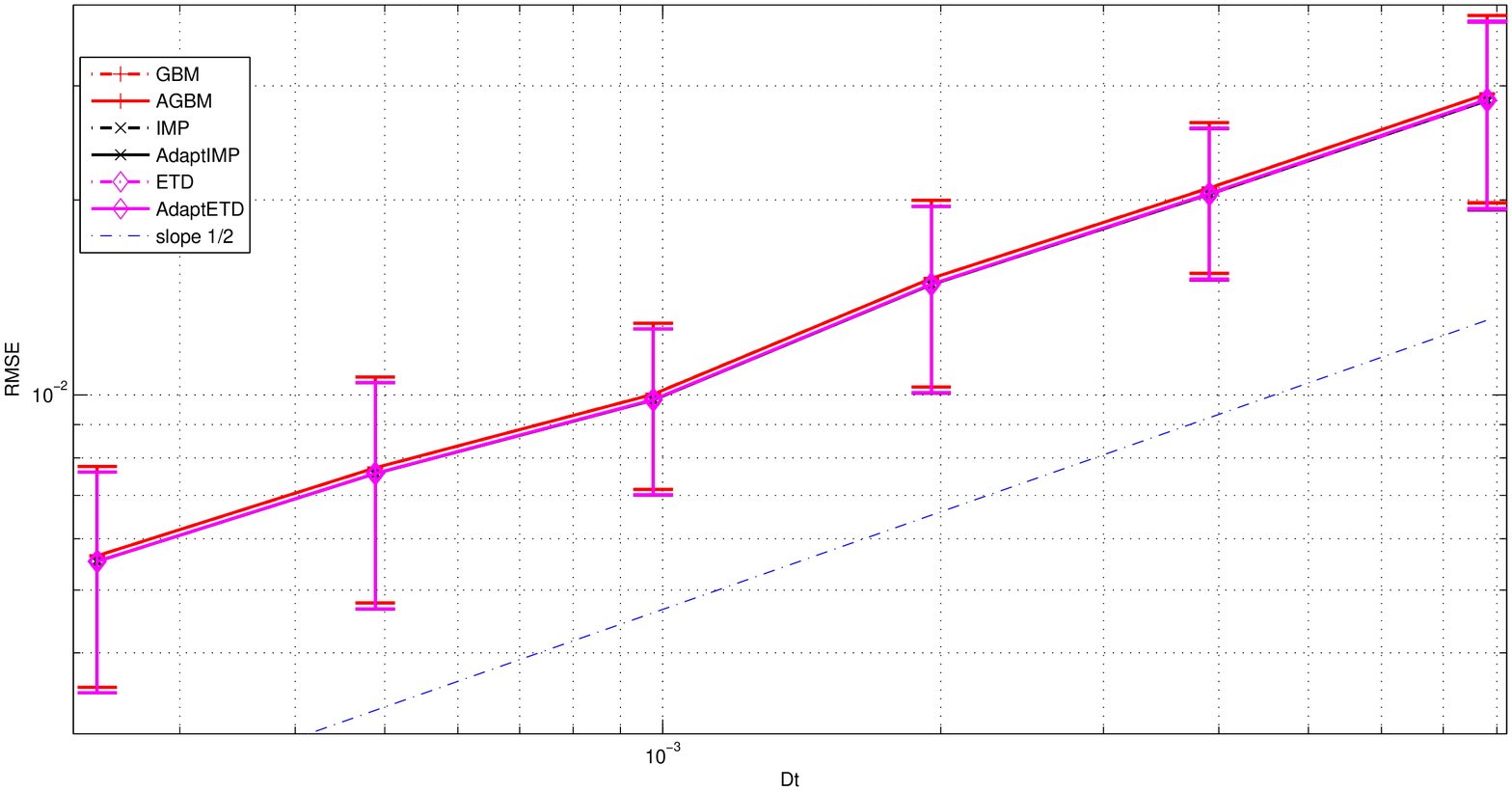}
\includegraphics[width=0.495\textwidth,height=0.45\textwidth]{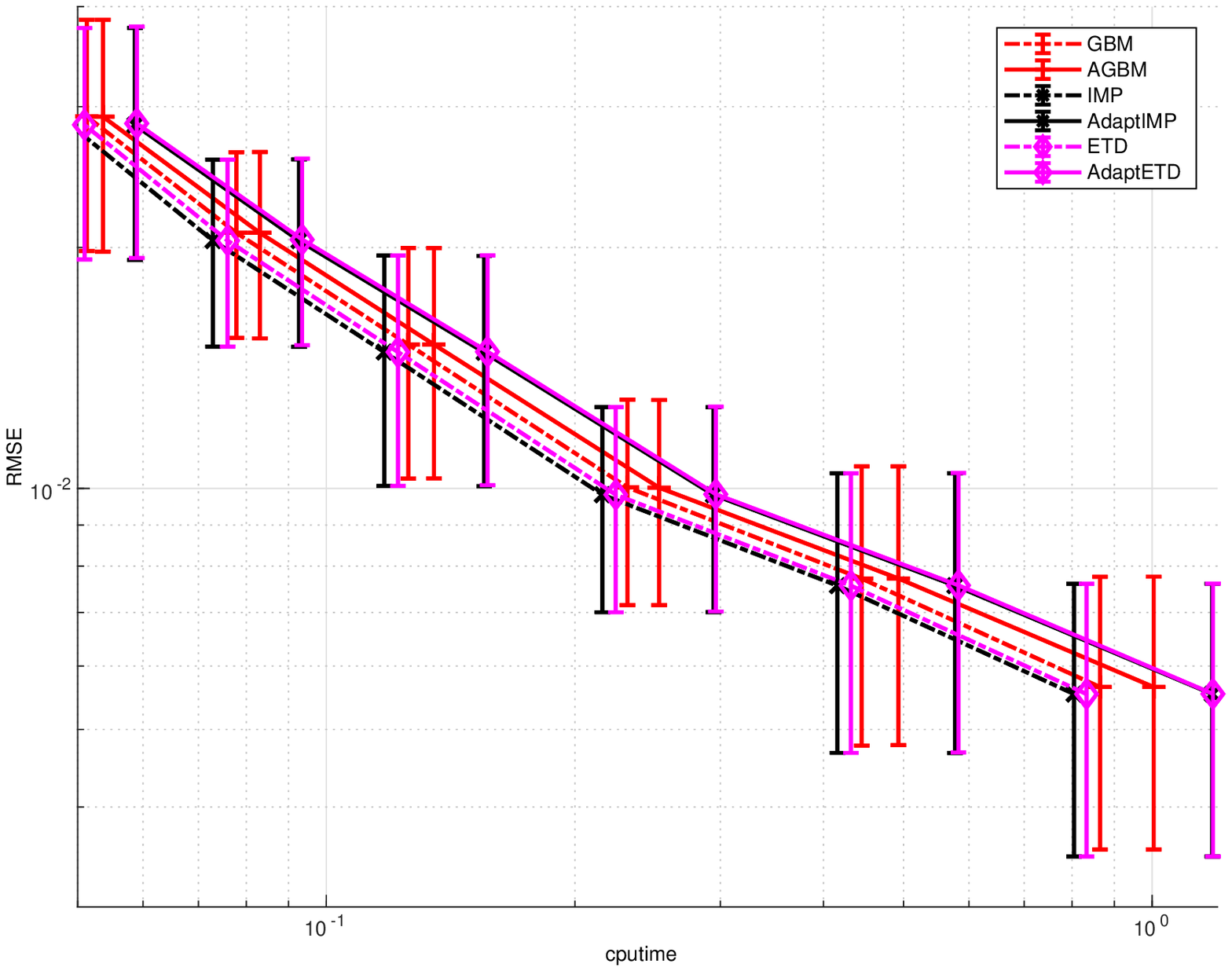}
\end{center}
\caption{SPDE \eqref{eq:SPDE} with  $\alpha=1$, $\beta=0.1$, $\gamma=1$ with $M=500$ samples (a) RMSE against $\Dt$ with a reference line with slope $\frac{1}{2}$.  (b)  RMSE against cputime.} 
  \label{fig:alpha1beta01}
\end{figure}

We now include nonlinearity in the noise and in Fig. \ref{fig:alpha05beta1} we take $\alpha=0.5$ and in Fig.  \ref{fig:alpha1beta01} we have $\alpha=1$. For a weaker nonlinearity in Fig. \ref{fig:alpha05beta1} we see the GBM methods are the most competitive.
 When the noise is strongly nonlinear in Fig. \ref{fig:alpha1beta01} we see all  the methods display the same convergence and in efficiency there is a slight advantage to the adaptivity.

In Fig. \ref{fig:alpha05beta1} although adaptive methods are more accurate, fixed step size methods are ahead in terms of CPU times. 

\appendix

\section{Proof of results from \secref{sec:boundedmoment}}
\label{app:boundedmoment}

\begin{proof}[Proof of Lemma \ref{lem:boundedness1}]
On $\Omega_{n+1}^N$ by construction we have that
$\norm{\sum_{i=1}^m B_i\Delta W_{t_n}^i}\leq 1$
for $n=0,1,\hdots, N-1$. Therefore, $\norm{ \SG{t_{n+1},t_n} }< \infty$ for  $n=0,1,\hdots,N-1$. 
We prove the Lemma on two  subsets of $\Omega_{n+1}^N$.
\begin{align}
1) & \quad S^{(1)}_{n+1}:=\Omega_{n+1}^N\cap \lbrace \omega \in \Omega \vert \norm{Y_n^N (\omega)} \leq 1 \rbrace\\    
2) & \quad S^{(2)}_{n+1}:=\Omega_{n+1}^N\cap \lbrace \omega \in \Omega \vert 1 \leq \norm {Y_n^N (\omega)} \leq  N^{1/(2c)} \rbrace.
\end{align}
First, on  $S^{(1)}_{n+1}$, we have from \eqref{eq:tamedEI0} and the triangle inequality that
\begin{equation}
\norm{Y_{n+1}^N} \leq \norm{ \SG{t_{n+1},t_n} } \left( \norm{Y_n^N}+\Delta t \norm{\tilde{F}(Y_n^N)}\right).
\end{equation}
Since $\norm{Y_n^N} \leq 1$   on $S^{(1)}_{n+1}$, and by the taming inequality \eqref{eq:tamebound}, we have that 
\begin{equation*}
\norm{Y_{n+1}^N} \leq \norm{ \SG{t_{n+1},t_n} } \left( 1+\Delta t \norm{F(Y_n^N)}\right).
\end{equation*}
Adding and subtracting $F(0)$ and applying the triangle inequality we get 
\begin{equation}
\label{eq:tmp0}
\norm{Y_{n+1}^N} \leq \norm{ \SG{t_{n+1},t_n} } \left( 1+\Delta t \norm{F(Y_n^N)-F(0)}+\Delta t \norm{F(0)}\right).
\end{equation}
Applying the Mean Value Theorem for $F$ and using the growth bound on the gradient of $F$ from Assumption \ref{ass:1} 
\begin{align}
 &\norm{F(Y_n^N)-F(0)}=\norm{\vec{D}F(\theta Y_n^N) {Y}_n^N}  \leq K \left( 1+\norm{Y_n^N} ^c \right) \norm{Y_n^N },
 \end{align}
where $0<\theta<1$. Thus, from \eqref{eq:tmp0} and due to the fact that $\Delta t< T$, we have 
\begin{equation*}
 \norm{Y_{n+1}^N}\leq \norm{\SG{t_{n+1},t_n}} \left( 1+T  K   \left( 1+\norm{Y_n^N} ^c \right) \norm{Y_n^N}  +T \norm{F(0)}\right). 
\end{equation*}
Now using that $\norm{Y_n^N} \leq 1$  and $ \norm{\sum_{i=1}^m B_i\Delta W_{t_k}^i}\leq 1$, we get on $S^{(1)}_{n+1}$
 \begin{equation} \label{eq:Y_bound_S1}
  \norm{Y_{n+1}^N} \leq \norm{\SG{t_{n+1},t_n}} \left( 1+2 T K+T \norm{F(0)}\right) \leq  \lambda.
 \end{equation}
Secondly : on the set $S^{(2)}_{n+1}$, we start from \eqref{eq:tamedEI0} by squaring the norm,
\begin{align}
   \norm{Y_{n+1}^N}^2 &\leq  \norm{ \SG{t_{n+1},t_n} }^2 \norm{Y_n^N+\Delta t \tilde{F}(Y_n^N)}^2 \nonumber \\
   &\leq \norm{ \SG{t_{n+1},t_n} }^2 \left( \norm{Y_n^N}^2+2 \Delta t^2 
\norm{ {F}(Y_n^N)}^2 +2 \Delta t  \innerproduct{Y_n^N}{ F(Y_n^N)} \right)
\label{eq:tmplem1}
\end{align}
where we have again used \eqref{eq:tamebound} on the last two terms. 
By the polynomial growth bound on $\vec{D}F$ (see Assumption \ref{ass:1})  on  $S^{(2)}_{n+1}$,($1 \leq \norm{Y_n^N (\omega)}\leq N^{(1/(2c)}$)  we have
 \begin{align*}
 \norm{F(Y_n^N)}^2& \leq  \left( \norm{F(Y_n^N)-F(0)}+\norm{F(0)}\right)^2 \\
 & \leq N   \left(2K+\norm{F (0)} \right)^2  \norm{Y_n^N} ^{2}.
 \end{align*}
The one sided Lipschitz condition on $F$ (Assumption \ref{ass:1}) and  Cauchy–Schwarz inequality give that
\begin{align*}
\innerproduct{Y_n^N}{F(Y_n^N)} & \leq \innerproduct{Y_n^N}{F(Y_n^N)-F(0)} +\innerproduct{Y_n^N}{F(0)}\\
& \leq  \left( K + \norm{F (0)}\right) \norm{Y_n^N} ^2,
\end{align*} 
where we have used that $1\leq \|Y_n^N\|$.
As a result, and since  $\Delta t =T/N$, we see that \eqref{eq:tmplem1} becomes
\begin{align} \label{eq:Y_bound_S2}
\begin{split}
 \norm{Y_{n+1}^N}^2 & \leq \norm{\SG{t_{n+1},t_n}}^2  \norm{Y_n ^N}^2 \left(1+\frac{2}{N} \left((2TK+T \norm{F (0)})^2+(TK+T \norm{F (0)})\right)  \right) \\
 & \leq \norm{\SG{t_{n+1},t_n}}^2  \norm{Y_n^N}^2(1+2\lambda/N)
  \leq \norm{\SG{t_{n+1},t_n}}^2  \norm{Y_n^N}^2 e^{2\lambda /N}.
\end{split}
\end{align}
The base  case of induction for  $n=0$  is obvious by initial condition on $\Omega_0^N=\Omega$. Let  $l \in \lbrace 0,1,\hdots N-1 \rbrace$ and assume $\norm{Y_{n}(\omega)}\leq D_n ^N(\omega)$  holds for all $n \in {0,1,\hdots,l}$ where  $\omega \in \Omega_n ^N$. We now  prove  that   
 $$\norm{Y_{l+1}(\omega)}\leq D_{l+1} ^N(\omega)$$
 for all $\omega \in \Omega_{l+1} ^N$.  For  all   $ \omega \in \Omega_{l+1} ^N$, we have $\norm{Y_{n}(\omega)}\leq D_n ^N(\omega)\leq N^{1/(2c)}$ ,   $n \in {0,1,\hdots,l}$  by induction hypothesis and $\Omega_{l+1} ^N \subseteq  \Omega_{n+1} ^N$ .  
For any $ \omega \in \Omega_{l+1} ^N$,  $\omega$ belongs to $S_{n+1} ^{(1)}$ or $S_{n+1} ^{(2)}$. For inductive argument we define  a random variable 
$$
\tau_l ^N (\omega):=\max \left(  \lbrace  -1 \rbrace \cup \lbrace  n \in \lbrace 0 ,1,\hdots,l-1 \rbrace \bigg\vert  \norm{Y_n (\omega)} \leq 1 \rbrace  \right)
$$ 
as done in \citep{hutzenthaler2012}.
This definition implies that $1 \leq \norm{Y_n(\omega)} \leq N^ {1/(2c)}$ for  all  $n \in \lbrace \tau_{l+1} ^N (\omega)+1, \tau_{l+2} ^N (\omega),\hdots, l \rbrace$. By estimation \ref{eq:Y_bound_S2},
\begin{align*}
  \norm{Y_{l+1} ^N  (\omega)} & \leq  \norm{\SG{t_{l+1},t_l} (\omega)} \norm{Y_{l} ^N (\omega)}    e^{\lambda /N}\\
&   \leq \hdots \leq  \norm{Y_{ \tau_{l+1} ^N (\omega)+1} ^N (\omega)}   \prod_{n=\tau_{l+1} ^N (\omega)+1 } ^l \left( \norm{\SG{t_{n+1}},t_{n} (\omega)  }  e^{\lambda /N}  \right)\\
 & \leq  \norm{Y_{ \tau_{l+1} ^N (\omega)+1} ^N (\omega)}   e^{\lambda } \sup _{u \in \lbrace 0,1,\hdots,l+1 \rbrace } \prod_{n=u } ^l \norm{\SG{t_{n+1}},t_{n} (\omega)  }   
\end{align*}
By considering \eqref{eq:Y_bound_S1}, the following completes  the induction step and proof. 
$$
\norm{Y_{l+1} ^N  (\omega)} \leq (\lambda + \norm{\xi (\omega)} )  e^{\lambda } \sup _{u \in \lbrace 0,1,\hdots,l+1 \rbrace } \prod_{n=u } ^l \norm{\SG{t_{n+1}},t_{n} (\omega)  }  
$$
\end{proof}

\begin{proof}[Proof of Lemma \ref{lem:Dboundedness}]
By H\"{o}lder's inequality and independence of the  Brownian increments
\begin{align*}
 \sup_{N \in \mathbb{N} }   \LpnormR{\sup_{n \in \{ 0,1,\hdots N \} }  D_n^N }& \leq  e^{\lambda}\left( \lambda + \Ltwopnorm{\xi} \right) \left(   \sup_{N \in \mathbb{N} }  \LtwopnormR{ \prod_{k=0}^{N-1}\norm{\SG {t_{k+1},t_k}}} \right) \nonumber \\ 
 & =  e^{\lambda}\left( \lambda + \Ltwopnorm{\xi} \right) \left(   \sup_{N \in \mathbb{N} } \prod_{k=0}^{N-1} \eval{\norm{ \SG {t_{k+1},t_k} }^{2p}}^{\frac{1}{2p}} \right).
\end{align*}
By Lemma \ref{lem:SGBoundedness}, we obtain the result.
\end{proof}

\begin{proof}[Proof of Lemma \ref{lem:complementProbability}]
For all $N \in \mathbb{N}$ and $q \in [1,\infty)$,  we have
\begin{align*}
\mathbb{P} [(\Omega _N^N)^c]& \leq  \mathbb{P} [  \sup _{0\leq k \leq N-1} D_k^N > N^{1/2c} ] + N  \mathbb{P} [  \norm{\sum_{i=1}^m B_i W_{t_1}^i} > 1    ] \\
& \leq   \eval{  \sup _{0\leq k \leq N-1} \vert D_k^N \vert  ^q } N^{-q/(2c)} + N  \mathbb{P} [  \norm{\sum_{i=1}^m B_i W_{T}^i} > \sqrt{N}    ] \\
&\leq   \eval{  \sup _{0\leq k \leq N-1} \vert D_k^N \vert  ^q } N^{-q/(2c)} +   m^{p-1} \max_{i \in \lbrace1,2, \hdots m\rbrace} \norm{B_i}^q  m \eval {  \norm{ W_{T}^1}^1     } N^{1-q/2}
\end{align*}
by definition of $\Omega _n^N $ in \eqref{eq:Omega}  and subadditivity  of probability measure and Markov inequality. The bounded moments of $D_n^N$ is given in \eqref{lem:Dboundedness}. For  given $p$,  $q $ is chosen such that  $(N^p \mathbb{P} [(\Omega _N^N)^c] )$ becomes bounded for any $N \in \mathbb{N}$
\end{proof}

\begin{proof}[Proof of Theorem \ref{teo:boundedMoment}]
The iterated numerical solution is given by
\begin{equation} \label{eq:sum_discrete_int}
Y_n^N= \SG{t_n,0} \xi 
+ \sum_{k=0}^ {n-1}  \SG{t_n,t_k}   \tilde{F} (Y_k^N) \Delta t.
\end{equation}
By taking the norm  in  the  space $L^p(\Omega,\mathbb{R}^d)$, applying the triangle inequality and noticing $ \Lpnorm{\tilde{F} (Y_k^N) \Delta t }\leq 1$, we  have
\begin{align} 
\Lpnorm{Y_n^N}& \leq \Lpnorm{\SG{t_n,0} \xi} 
+ \Lpnorm{\sum_{k=0}^ {n-1}  \SG{t_n,t_k}   \tilde{F} (Y_k^N) \Delta t} \leq  C (\norm{\xi}+N) 
\label{eq:tmp1T1}
\end{align}
where $C=  \sup_{ s,t \in [0,T]}   \Lpnorm{\SG{s,t}}$.
The  existence of $N$ on the RHS of the inequality \eqref{eq:tmp1T1} appears to be an obstacle in completing the proof. However, following the techniques of \  \cite[equation (59)]{hutzenthaler2012}]  
\begin{align}
\sup_{N \in \mathbb{N}} \sup_{0\leq n\leq N}  \Lpnorm{1_{(\Omega_n^N)^c} Y_n^N}  & < \infty
 \label{eq:tmpPt1}
\end{align}
and using Lemmas \ref{lem:boundedness1} and \ref{lem:Dboundedness}
  \begin{equation}
  \sup_{N \in \mathbb{N}} \sup_{0\leq n \leq N}  \Lpnorm{1_{(\Omega_n^N)} Y_n^N } \leq   \sup_{N \in \mathbb{N}} \sup_{0\leq n\leq N} \Lpnorm{D_n^N } < \infty.
  \label{eq:tmpPt2}
  \end{equation}
  Combining \eqref{eq:tmpPt1} and \eqref{eq:tmpPt2} completes the proof.
\end{proof}

\section*{Acknowledgements}
\noindent
We would like to thank Nesin Mathematics Village, Nesin Mathematik K\"oy\"u, for hosting us and support that allowed this work to take place. The research of the first author was supported by grant number 120F299 under TUBITAK 1002.

\bibliographystyle{unsrtnat}
\bibliography{taming_references}  

\end{document}